\documentclass[12pt]{article}

\usepackage{epstopdf}
\usepackage{amsmath,amssymb,bm} 
\usepackage{amsthm}
\usepackage{amsmath,amssymb,amsfonts,amsthm,bm,mathrsfs}

\newtheorem{thm}{Theorem}[section]

\newtheorem{lem}[thm]{Lemma}

\newtheorem{con}[thm]{Conjecture}

\newtheorem{pf}{proof}[section]
\theoremstyle{remark}

\makeatletter \@addtoreset{equation}{section} \makeatother
\makeindex 
\def\qed{\hfill \rule{4pt}{7pt}}

\begin{document}
\begin{center}

 {\Large \bf Higher Order Tur\'an Inequalities for the\\[5pt] Distinct Partition Function}

\end{center}

\begin{center}
Janet J.W. Dong$^{1}$ and {Kathy Q. Ji}$^{2}$   \vskip 2mm

$^{1,\,2}$ Center for Applied Mathematics\\[3pt]
Tianjin University\\[3pt]
Tianjin 300072, P.R. China\\[6pt]
   \vskip 2mm

 Emails: $^1$dongjinwei@tju.edu.cn and   $^2$kathyji@tju.edu.cn
\end{center}

\vskip 6mm \noindent {\bf Abstract.}
We prove that the number  $q(n)$ of partitions into distinct parts  is log-concave for   $n \geq 33$ and satisfies   the higher order Tur\'an inequalities for $n\geq 121$ conjectured by  Craig and Pun.  In doing so,  we establish explicit error terms   for $q(n)$ and for $q(n-1)q(n+1)/q(n)^2$ based on Chern's asymptotic formulas  for $\eta$-quotients.



 \vskip 6mm

\section{Introduction}

The objective of this paper is to explore the Tur\'an inequalities and the higher order Tur\'an inequalities for the distinct partition  function.  The Tur\'an inequalities (also called Newton inequalities) and the higher order Tur\'an inequalities (also called the cubic inequalities) arise in the study of the Maclaurin coefficients of real entire functions in Laguerre-P\'olya class, see, for example,  \cite{Dimitrov1998} and \cite{Szego-1948}. A sequence $\{\alpha_n\}_{n\geq 0}$ of real numbers is said to satisfy the    Tur\'an inequalities or to be log-concave if for $n\geq 1$,
\begin{equation}\label{log-con}
  \alpha^2_n\geq \alpha_{n-1}\alpha_{n+1}.
\end{equation}
The sequence $\{\alpha_n\}_{n\geq 0}$ is said to satisfy the higher order Tur\'an inequalities if for $n \geq 1$,
\begin{equation}\label{hig-Turan}
 4(\alpha^2_n-\alpha_{n-1}\alpha_{n+1})(\alpha^2_{n+1}-\alpha_n\alpha_{n+2})\geq (\alpha_{n}\alpha_{n+1}-\alpha_{n-1}\alpha_{n+2})^2.
\end{equation}
The    Tur\'an inequalities and the higher order Tur\'an inequalities are closely related to   the Jensen polynomials, see, for example,  \cite{Craven-Csordas-1989, Csordas-Norfolk-Varga-1986,  Csordas-Varga-1990}.
  The  Jensen polynomials $J_{\alpha}^{d,n}(X)$ of degree $d$ and shift $n$ associated to the sequence $\{\alpha_n\}_{n\geq0}$ are defined by
\[
	J_{\alpha}^{d,n}(X)=\sum_{i=0}^{d} \binom{d}{i} \alpha_{n+i} X^i.
\]
For $d=2$ and shift $n-1$, the Jensen polynomial $J_{\alpha}^{2,n-1}(X)$ reduces to
\[J_{\alpha}^{2,n-1}(X)=\alpha_{n-1}+2\alpha_{n} X+\alpha_{n+1} X^2.
\]
It is clear that $\{\alpha_n\}_{n\geq0}$ is log-concave at $n$ if and only if  $J_{\alpha}^{2,n-1}(X)$ has only real roots. In general, we say that the sequence $\{\alpha_n\}_{n\geq0}$ satisfies the order $d$ Tur\'an inequality at $n$ if and only if $J_{\alpha}^{d,n-1}(X)$ is hyperbolic. Recall that a polynomial is hyperbolic if all of its roots are real.

 The  Tur\'an inequalities and the higher order Tur\'an inequalities for the partition function were initially investigated by   Chen \cite{Chen-2017}, DeSalvo and Pak \cite{DeSalvo-Pak-2015} and  Nicolas \cite{Nicolas-1978}. Recall that a partition of $n$ is a finite list of nondecreasing positive integers $\lambda=\left(\lambda_1, \lambda_2, \ldots, \lambda_r\right)$ such that $\lambda_1+\lambda_2+\cdots+\lambda_r=n$. Let $p(n)$ denote the number of partitions of $n$.  DeSalvo and Pak \cite{DeSalvo-Pak-2015} and  Nicolas \cite{Nicolas-1978} proved that the partition function $p(n)$ is log-concave
for $n \geq 26$. Chen \cite{Chen-2017}  conjectured that $p(n)$ satisfies the higher order Tur\'an inequalities
for $n \geq 95$, which was proved by Chen, Jia, and Wang \cite{Chen-Jia-Wang-2019}.  Chen, Jia, and Wang \cite{Chen-Jia-Wang-2019} further conjectured that for $d\geq 4$, there exists a
positive integer $N_p(d)$ such that $p(n)$ satisfies the order $d$  Tur\'an inequalities  for $n\geq N_p(d)$, that is,   the Jensen polynomial $J_{p}^{d,n-1}(X)$ associated to $p(n)$  is hyperbolic for $n\geq N_p(d)$. Griffin, Ono, Rolen, and Zagier \cite{Griffin-Ono-Rolen-Zagier-2019} confirmed their conjecture for   sufficiently large $n$.  More recently,   Tur\'an inequalities for  other partition functions have been extensively investigated, see, for example, Bringmann, Kane, Rolen, and Tripp \cite{Bringmann-Kane-Rolen-Tripp-2021}, Dong, Ji, and Jia \cite{Dong-Ji-Jia-2020}, Engel \cite{Engel-2017}, Liu and Zhang \cite{Liu-Zhang-2021}, Jia \cite{Jia-2022} and  Ono, Pujahari, and Rolen \cite{Ono-Pujahari-Rolen-2019}.

The   goal of this paper is to investigate the  Tur\'an inequalities and the higher order Tur\'an inequalities for the  distinct partition  function.  Let $q(n)$ denote the number of partitions of $n$ into distinct parts.
	For example, there are eight partitions of $9$ into distinct parts:
	\[(9),\,(8,1),\,(7,2),\,(6,3),\,(6,2,1),\,(5,3),\,(5,2,1),\,(4,3,2).\]
It's  known from Euler's partition theorem   that $q(n)$ also counts the number of partitions of $n$ into odd parts, see Andrews \cite[Chapter 1]{Andrews-1998} or Euler\cite{Euler-1748} .

	The generating function for $q(n)$ is given by
\begin{align*}
	\sum_{n\ge 0}q(n)q^n
	&=\prod_{j=1}^{\infty}(1+q^j)=\prod_{j=0}^{\infty}\frac{1}{1-q^{2j+1}}.
\end{align*}
Hagis \cite{ Hagis-1963} and Hua \cite{Hua-1942} established a Rademacher-type formula for $q(n)$ by the circle method in terms of Kloosterman sums and Bessel functions. Based on this formula,     Craig and Pun\cite{Craig-Pun-2021} showed that $q(n)$ satisfies the order $d$ Tur\'an inequalities  for sufficiently large $n$ by using   a general result of Griffin, Ono, Rolen, and Zagier \cite{Griffin-Ono-Rolen-Zagier-2019}. They also made the following  conjecture.
\begin{con}[Craig-Pun] \label{cp-1} The function  $q(n)$ is log-concave for $n\geq 33$ and satisfies the higher order Tur\'an inequalities for $n\geq 121$.
\end{con}

The main contribution  of this paper is to confirm
Conjecture \ref{cp-1}. Instead of using the  Rademacher-type formula for $q(n)$  due to Hagis \cite{ Hagis-1963} and Hua \cite{Hua-1942}, we find Chern's asymptotic formulas  for $\eta$-quotients \cite{Chern-2019} play an important role in the proof of Conjecture \ref{cp-1}.

	Appealing to Chern's asymptotic formulas, 	we obtain an asymptotic formula for $q(n)$ with an effective bound on the error term. To state our bound, we adopt the following notation:
	\begin{equation}\label{defi-x}
	    {\nu}(n)=\frac{\pi\sqrt{24n+1}}{6\sqrt{2}}.
	\end{equation}
	We have the following asymptotic formula for $q(n)$.
		\begin{thm}\label{eq-q-thm}For  ${\nu}(n)\geq 21$, or equivalently, $n\geq 135$,
		\begin{equation}\label{eq-q}
			q(n)=\frac{   \sqrt{2}\pi^2}{12 \nu(n)}I_{1}(\nu(n))+R(n),
		\end{equation}
		where $I_1(s)$ is the first modified Bessel function of the first kind defined as
\begin{equation}\label{eq-integral-bessel-1}
	I_1(s)=\frac{s}{\pi}\int_{-1}^1(1-t^2)^{\frac{1}{2}}e^{st}{\rm{d}}t,
\end{equation}
		and
		\begin{equation}\label{bound-R}
				|R(n)|\leq \frac{\sqrt{3} \pi^{\frac{3}{2}}}{  6\nu(n)^{\frac{1}{2}}}  \exp\left(\frac{\nu(n)}{3}\right) .
		\end{equation}
	\end{thm}
Using Theorem \ref{eq-q-thm}, we establish an upper bound and a lower bound for $q(n)$ which are required  in the proof of Conjecture \ref{cp-1}.
	
\begin{thm}\label{lem-M-1}   Let
 \begin{equation}\label{defi-M}
 	M(n):=\frac{\sqrt{2}\pi^2}{12\nu(n)}I_1(\nu(n)).
 \end{equation}  For $\nu(n) \geq 43 $,
	\begin{equation}\label{eq-main}
		M(n)\left(1-\frac{1}{\nu(n)^6}\right)\leq q(n)\leq M(n)\left(1+\frac{1}{\nu(n)^6}\right).
	\end{equation}
	
\end{thm}
	
Let
\begin{equation}\label{defi-Bn}
	Q(n)=\frac{q(n-1)q(n+1)}{q(n)^2}.
     \end{equation}
It's evident from \eqref{log-con} that $q(n)$ is log-concave for $n\geq 33$ is equivalent to $Q(n)\leq 1$ for $n\geq 33$. Using \eqref{hig-Turan}, one can check that   $q(n)$ satisfies  the higher order Tur\'an inequalities for $n\geq121$ whenever
\begin{equation}\label{eq-q-turan}
4(1-Q(n))(1-Q(n+1))-(1-Q(n)Q(n+1))^2\geq0
	\end{equation}
 for $n\geq 121$. Hence to prove Conjecture \ref{cp-1}, we first establish  an efficient lower and upper bounds for $Q(n)$ by using Theorem \ref{lem-M-1}.
\begin{thm}\label{thm-B}  Let
\begin{equation}\label{defi-I-e}
    E_Q(n):=1-\frac{\pi ^4}{36 \nu(n)^3}+\frac{\pi ^4}{12 \nu(n)^4}-\frac{\pi ^4}{32 \nu(n)^5}.
\end{equation}
Then for $\nu(n)\geq 67$,
	\begin{equation}\label{thm-B-low}
		E_Q(n)-\frac{135 }{\nu(n)^6}<Q(n)<E_Q(n)+\frac{126+\frac{\pi ^8}{1296} }{\nu(n)^6}.
	\end{equation}
\end{thm}
It turns out that  Conjecture \ref{cp-1} follows from
 Theorem \ref{thm-B}. We therefore arrive at the following consequence.
\begin{thm}\label{coro-tulan-equiv}
	 For $n\geq121$, the cubic polynomial
	 \[q(n-1)+3q(n)x+3q(n+1)x^2+q(n+2)x^3\]
	 has  only  real zeros.
\end{thm}

The paper is organized as follows.   In Section 2, we   prove some inequalities involving the first modified Bessel function of the first kind, which are necessary in the proof of Theorem \ref{thm-B}. In Section 3, we prove Theorem \ref{eq-q-thm} with the aid of Chern's asymptotic formulas  for $\eta$-quotients and then use  Theorem \ref{eq-q-thm} to prove Theorem \ref{lem-M-1}. Section 4 is devoted to the proof of Theorem \ref{thm-B} with the aid of Theorem \ref{eq-q-thm} and the inequalities on the first modified Bessel function of the first kind established in Section 2.     In Section 5,  we give a proof of Conjecture \ref{cp-1} by using Theorem \ref{thm-B}. We  conclude in Section 6 with some problems for further investigation.

\section{Explicit bounds for $I_1(s)$}

In order to make the asymptotic formula in Theorem \ref{defi-x} useful in the proof of  Conjecture \ref{cp-1},  we need some inequalities on the first modified Bessel function $I_1(s)$  of the first kind. Before doing this, let us first recall the definitions of the Gamma function $\Gamma(a)$ and the upper incomplete Gamma function  $\Gamma(a,s)$, see \cite[Chapter 6]{Abramowitz-Stegun-1972}.

The Gamma function $\Gamma(a)$ is   defined by
\[\Gamma(a)=\int_{0}^{\infty}t^{a-1}e^{-t}{\rm{d}}t.\]
It's known that
\[\Gamma(a+1)=a\Gamma(a),\quad \Gamma\left(\frac{1}{2}\right)=\sqrt{\pi},\]
see \cite[p. 32--34]{Rademacher-1973}.

The upper incomplete Gamma function  $\Gamma(a,s)$ is defined by
\[\Gamma(a,s)=\int_{s}^{\infty}t^{a-1}e^{-t}{\rm{d}}t.\]
The following estimate on $\Gamma(a,s)$ can be derived from   the proof of Proposition $2.6$ of Pinelis \cite{Pinelis-2020} which is required in the proof of Lemma \ref{thm-bessel-1}. For $a\ge1$ and $s\ge a$,
\begin{equation}\label{p-incomgamma}
    \Gamma(a,s)\leq as^{a-1}e^{-s}.
\end{equation}

The first inequality on $I_1(s)$ for the proof of  Conjecture \ref{cp-1} is due to  Bringmann, Kane, Rolen, and Trippin \cite[Lemma 2.2]{Bringmann-Kane-Rolen-Tripp-2021}.

\begin{lem}[Bringmann-Kane-Rolen-Trippin] For $s\geq 1$,
\begin{equation}\label{ineq-BKRT}
	I_1(s) \leq \sqrt{\frac{2}{\pi s}} e^s.
\end{equation}
\end{lem}
We also need the further estimate   on $I_1(s)$.

\begin{lem}\label{thm-bessel-1}
Let
\begin{equation}\label{defi-I-e}
    E_I(s):=1-\frac{3}{8s}-\frac{15}{128s^2}-\frac{105}{1024s^3}
		-\frac{4725}{32768s^4}-\frac{72765}{262144s^5}.
\end{equation}
Then for $s\geq 26$,
	\begin{align}
	\frac{e^s}{\sqrt{2\pi s}}\left( E_I(s)-\frac{31}{s^6}\right)	\leq	I_1(s) \leq \frac{e^s}{\sqrt{2\pi s}}\left( E_I(s)+\frac{31}{s^6}\right), \label{ineqBessel-1}
\end{align}
\end{lem}
{\noindent \it Proof.}
We start with the integral definition \eqref{eq-integral-bessel-1} of $I_1(s)$,
\begin{equation}\label{eq-integral-bessel-2}
	I_1(s)=\frac{s}{\pi}\int_{0}^1(1-t^2)^{\frac{1}{2}}e^{st}{\rm{d}}t+\frac{s}{\pi}\int_{-1}^0(1-t^2)^{\frac{1}{2}}e^{st}{\rm{d}}t.
\end{equation}
It is clear that
\begin{equation}\label{eq-integral-bessel-r1}
	\left|\frac{s}{\pi}\int_{-1}^0(1-t^2)^{\frac{1}{2}}e^{st}{\rm{d}}t\right|\leq \frac{s}{\pi}.
\end{equation}
We next estimate the first integral in \eqref{eq-integral-bessel-2}.
Setting $u=1-t$, we have
\begin{align}\label{eq-integral-bessel-3}
\frac{s}{\pi}\int_{0}^1(1-t^2)^{\frac{1}{2}}e^{st}{\rm{d}}t&=\frac{se^s}{\pi}\int_{0}^1(2-u)^{\frac{1}{2}}u^{\frac{1}{2}}e^{-su}{\rm{d}}u.
\end{align}
By Taylor's formula, we find that
\begin{align}\label{taylor}
(2-u)^{\frac{1}{2}}=\sqrt{2}-\frac{u}{2 \sqrt{2}}-\frac{u^2}{16 \sqrt{2}}-\frac{u^3}{64 \sqrt{2}}-\frac{5u^4}{1024 \sqrt{2}}-\frac{7u^5}{4096 \sqrt{2}}+c(\xi)u^6,
\end{align}
where
\begin{equation}\label{eq-c(u)}
	c(\xi)=\frac{1}{6!}\left(\frac{{\rm{d}}^6}{{\rm{d}}u^6}(2-u)^{\frac{1}{2}}\right)_{u=\xi}=-\frac{21}{1024}(2-\xi)^{-\frac{11}{2}} \quad 	\text{for some} \quad  \xi\in(0,1).
	\end{equation}

Substituting \eqref{taylor} into  \eqref{eq-integral-bessel-3}, we obtain
\begin{align}
&	\frac{\sqrt{2}se^s}{\pi}\int_{0}^1\left(u^{\frac{1}{2}}-\frac{1}{4}u^{\frac{3}{2}}-\frac{1}{32}u^\frac{5}{2}-\frac{1}{128}u^\frac{7}{2}-\frac{5}{2048}u^\frac{9}{2}-\frac{7}{8192}u^\frac{11}{2}\right)e^{-su}{\rm{d}}u\nonumber\\[6pt]
&\qquad+\frac{se^s}{\pi}\int_{0}^1c(\xi)u^{\frac{13}{2}}
	e^{-su}{\rm{d}}u\nonumber\\[6pt]
&=\frac{\sqrt{2}se^s}{\pi}\left(\int_{0}^{\infty}-\int_1^{\infty}\right)\left(u^{\frac{1}{2}}-\frac{1}{4}u^{\frac{3}{2}}-\frac{1}{32}u^\frac{5}{2}-\frac{1}{128}u^\frac{7}{2}-\frac{5}{2048}u^\frac{9}{2}-\frac{7}{8192}u^\frac{11}{2}\right)e^{-su}{\rm{d}}u\nonumber\\[6pt]
&\qquad+\frac{se^s}{\pi}\int_{0}^1c(\xi)u^{\frac{13}{2}}
e^{-su}{\rm{d}}u:=I^{(1)}_1(s)+I^{(2)}_1(s)+I^{(3)}_1(s).
\label{eq-integral-bessel-5}
\end{align}
Evaluating the first integral $I^{(1)}_1(s)$ in \eqref{eq-integral-bessel-5} yields the main term:
\begin{align}
I^{(1)}_1(s)=&\frac{\sqrt{2}se^s}{\pi}\int_{0}^{\infty}\left(u^{\frac{1}{2}}-\frac{1}{4}u^{\frac{3}{2}}-\frac{1}{32}u^\frac{5}{2}-\frac{1}{128}u^\frac{7}{2}-\frac{5}{2048}u^\frac{9}{2}-\frac{7}{8192}u^\frac{11}{2}\right)e^{-su}{\rm{d}}u\nonumber\\[6pt]
&=\frac{\sqrt{2}e^s}{\sqrt{s}\pi}\int_{0}^{\infty}\left((su)^{\frac{1}{2}}-\frac{1}{4s}(su)^{\frac{3}{2}}-\frac{1}{32s^2}(su)^\frac{5}{2}-\frac{1}{128s^3}(su)^\frac{7}{2}\right.\nonumber\\[6pt]
&\phantom{=\;\;}\left.
\qquad-\frac{5}{2048s^4}(su)^\frac{9}{2}-\frac{7}{8192s^{5}}(su)^\frac{11}{2}\right)e^{-su}{\rm{d}}(su)\nonumber\\[6pt]
&=\frac{\sqrt{2}e^s}{\sqrt{s}\pi}\left(\Gamma\left(\frac{3}{2}\right)-\frac{1}{4s}\Gamma\left(\frac{5}{2}\right)-\frac{1}{32s^2}\Gamma\left(\frac{7}{2}\right)-\frac{1}{128s^3}\Gamma\left(\frac{9}{2}\right)\right.\nonumber\\[6pt]
&\phantom{=\;\;}\left.
\qquad-\frac{5}{2048s^4}\Gamma\left(\frac{11}{2}\right)-\frac{7}{8192s^5}\Gamma\left(\frac{13}{2}\right)\right)\nonumber\\[6pt]
&=\frac{e^s}{\sqrt{2\pi s}}\left(1-\frac{3}{8s}-\frac{15}{128s^2}-\frac{105}{1024s^3}-\frac{4725}{32768s^4}-\frac{72765}{262144s^5}\right).
\label{eq-integral-bessel-main}
\end{align}

We proceed to evaluate the second integral $I^{(2)}_1(s)$ in \eqref{eq-integral-bessel-5}:
\begin{align}
I^{(2)}_1(s)=&-\frac{\sqrt{2}se^s}{\pi}\int_{1}^{\infty}\left(u^{\frac{1}{2}}-\frac{1}{4}u^{\frac{3}{2}}-\frac{1}{32}u^\frac{5}{2}-\frac{1}{128}u^\frac{7}{2}-\frac{5}{2048}u^\frac{9}{2}-\frac{7}{8192}u^\frac{11}{2}\right)e^{-su}{\rm{d}}u\nonumber\\[6pt]
	&=-\frac{\sqrt{2}e^s}{\sqrt{s}\pi}\int_{s}^{\infty}\left((su)^{\frac{1}{2}}-\frac{1}{4s}(su)^{\frac{3}{2}}-\frac{1}{32s^2}(su)^\frac{5}{2}-\frac{1}{128s^3}(su)^\frac{7}{2}-\right.\nonumber\\[6pt]
	&\phantom{=\;\;}\left.
	\qquad-\frac{5}{2048s^4}(su)^\frac{9}{2}-\frac{7}{8192s^{5}}(su)^\frac{11}{2}\right)e^{-su}{\rm{d}}(su)\nonumber\\[6pt]
	&=\frac{\sqrt{2}e^s}{\sqrt{s}\pi}\left(-\Gamma\left(\frac{3}{2},s\right)+\frac{1}{4s}\Gamma\left(\frac{5}{2},s\right)+\frac{1}{32s^2}\Gamma\left(\frac{7}{2},s\right)+\frac{1}{128s^3}\Gamma\left(\frac{9}{2},s\right)\right.\nonumber\\[6pt]
	&\phantom{=\;\;}\left.
	\qquad+\frac{5}{2048s^4}\Gamma\left(\frac{11}{2},s\right)+\frac{7}{8192s^5}\Gamma\left(\frac{13}{2},s\right)\right),\nonumber\\[6pt]
	&\overset{\eqref{p-incomgamma}}{\leq}  \frac{\sqrt{2}e^s}{\sqrt{s}\pi}\left(\frac{3}{2}+\frac{1}{4}\cdot\frac{5}{2}+\frac{1}{32}\cdot\frac{7}{2}+\frac{1}{128}\cdot\frac{9}{2}+\frac{5}{2048}\cdot\frac{11}{2}+\frac{7}{8192}\cdot\frac{13}{2} \right)\sqrt{s}e^{-s}\nonumber\\[6pt]
	&=\frac{37495\sqrt{2}}{16384\pi} \quad \quad  \text{for}    \quad  s\geq \frac{13}{2}.
	\label{eq-integral-bessel-r2}
\end{align}

It remains to   estimate $I^{(3)}_1(s)$ in \eqref{eq-integral-bessel-5}. From \eqref{eq-c(u)}, we see that
\begin{align}\label{eq-integral-bessel-r3}
|I^{(3)}_1(s)|=\left|\frac{se^s}{\pi}\int_{0}^1c(\xi)u^{\frac{13}{2}}e^{-su}{\rm{d}}u\right|&\leq\frac{21s^{-\frac{13}{2}}e^s}{1024\pi}\int_0^{\infty}(su)^{\frac{13}{2}}e^{-su}{\rm{d}}(su)\nonumber\\[6pt]
&=\frac{21s^{-\frac{13}{2}}e^s}{1024\pi}\Gamma\left(\frac{15}{2}\right)\nonumber\\[6pt]
&=\frac{2837835}{131072\sqrt{\pi}}s^{-\frac{13}{2}}e^s.
\end{align}
Combining \eqref{eq-integral-bessel-r1}, \eqref{eq-integral-bessel-main},   \eqref{eq-integral-bessel-r2} and \eqref{eq-integral-bessel-r3}, we derive that for $s\geq \frac{13}{2}$,
\begin{equation}\label{eq-I-main}
I_1(s)=\frac{e^s}{\sqrt{2\pi s}}\left(1-\frac{3}{8s}-\frac{15}{128s^2}-\frac{105}{1024s^3}-\frac{4725}{32768s^4}-\frac{72765}{262144s^5}\right)+r(s),
\end{equation}
where
\begin{align*}\label{eq-r(s)}
	|r(s)|&\leq\frac{s}{\pi}+\frac{37495\sqrt{2}}{16384\pi}+\frac{2837835}{131072\sqrt{\pi}}s^{-\frac{13}{2}}e^s\nonumber\\[6pt]
	&=\frac{e^s}{\sqrt{2\pi s}}\cdot\frac{1}{s^6}\left(\left(\frac{\sqrt{2}s}{\sqrt{\pi}}+\frac{37495   }{8192\sqrt{\pi}}\right)s^{\frac{13}{2}}e^{-s}+\frac{2837835\sqrt{2}}{131072}\right).
\end{align*}
To prove \eqref{ineqBessel-1}, it suffices to show that  for $s\geq 26$,
\begin{equation}\label{eq-I-r}
|r(s)|\leq\frac{e^s}{\sqrt{2\pi s}}\cdot\frac{31}{s^6}.
\end{equation}
Define
\[f(s):=\left(\frac{\sqrt{2}s}{\sqrt{\pi}}+\frac{37495   }{8192\sqrt{\pi}}\right)s^{\frac{13}{2}}e^{-s}+\frac{2837835\sqrt{2}}{131072}.\]
Observe that
\[f'(s)=s^{\frac{11}{2}}e^{-s}\left(-\frac{\sqrt{2}}{\sqrt{\pi}}s^2+\frac{5 \left(12288 \sqrt{2}-7499\right)}{8192 \sqrt{\pi }}s +\frac{487435}{16384 \sqrt{\pi }}\right).\]
Since $f'(s)<0$ for $s\geq 8$, we deduce that $f(s)$ is decreasing when $s\geq 8$. This implies that when  $s\geq 26$,
\[f(s)\leq f(26)\approx 30.8068<31.\]
 Hence the inequality \eqref{eq-I-r} is valid.
 Combining \eqref{eq-I-main} and \eqref{eq-I-r}, we are led to \eqref{ineqBessel-1}   in Lemma \ref{thm-bessel-1}. This completes the proof.  \qed

Based on Lemma \ref{thm-bessel-1}, we obtain the following inequalities on $I_1(\nu(n-1))I_1(\nu(n+1))/I_1(\nu(n))^2$, which are required in the proof of Theorem \ref{thm-B}.
\begin{lem}\label{lem-A-1}
 For $\nu(n)\geq60$,

 \begin{align}\label{eq-lem-B1}
 	\frac{I_1(\nu(n-1))I_1(\nu(n+1))}{I_1(\nu(n))^2}
 	&\geq \frac{\nu(n)}{\sqrt{\nu(n-1)\nu(n+1)}}\left(1-\frac{ \pi^4}{36 \nu(n)^3}-\frac{5  \pi^8}{2592 \nu(n)^7}\right)\nonumber\\[6pt]
 	&~~~~~~~\times\left(1-\frac{ \pi^4}{32  \nu(n)^5}-\frac{129}{\nu(n)^6}\right),
 \end{align}
 \begin{align}\label{eq-lem-B2}
 	\frac{I_1(\nu(n-1))I_1(\nu(n+1))}{I_1(\nu(n))^2}
 	&\leq \frac{\nu(n)}{\sqrt{\nu(n-1)\nu(n+1)}}\left(1-\frac{ \pi^4}{36 \nu(n)^3}+\frac{ \pi^8}{1296 \nu(n)^6}\right) \nonumber\\[6pt]
 	&~~~~~~~\times\left(1-\frac{ \pi^4}{32  \nu(n)^5}+\frac{121}{\nu(n)^6}\right).
 \end{align}

\end{lem}
\proof
Using \eqref{ineqBessel-1}, we find that   for $\nu(n)\geq 26$,
\begin{align}\label{eq-I-1}
	\frac{I_1(\nu(n-1))I_1(\nu(n+1))}{I_1(\nu(n))^2}&\geq  \frac{\nu(n)}{\sqrt{\nu(n-1)\nu(n+1)}}e^{  \nu(n-1)+\nu(n+1)-2\nu(n)}L(n)
\end{align}
and
\begin{align}\label{eq-I-2}
	\frac{I_1(\nu(n-1))I_1(\nu(n+1))}{I_1(\nu(n))^2}&\leq  \frac{\nu(n)}{\sqrt{\nu(n-1)\nu(n+1)}}e^{  \nu(n-1)+\nu(n+1)-2\nu(n)}R(n),
\end{align}
where $\nu(n)$ and $E_{I}(s)$ are defined as \eqref{defi-x} and \eqref{defi-I-e} respectively,
\begin{align}\label{defi-tl}
L(n)=\frac{\left(E_{I}(\nu(n-1))-\frac{31}{\nu(n-1)^6}\right)\left(E_{I}(\nu(n+1))-\frac{31}{\nu(n+1)^6}\right)}{\left(E_{I}(\nu(n))+\frac{31}{\nu(n)^6}\right)^2}
\end{align}
and
\begin{align}\label{defi-hl}
R(n)=\frac{\left(E_{I}(\nu(n-1))+\frac{31}{\nu(n-1)^6}\right)\left(E_{I}(\nu(n+1))+\frac{31}{\nu(n+1)^6}\right)}{\left(E_{I}(\nu(n))-\frac{31}{\nu(n)^6}\right)^2}.
\end{align}

To obtain \eqref{eq-lem-B1} and \eqref{eq-lem-B2}, we intend to estimate $\exp({  \nu(n-1)+\nu(n+1)-2\nu(n)})$, $L(n)$ and $R(n)$ in terms of $\nu(n)$. From the definition \eqref{defi-x} of $\nu(n)$,  we see that for  $\nu(n)\geq 2$,
\begin{equation}\label{eq-x-Rel}
    \nu(n-1)=\sqrt{\nu(n)^2-\frac{\pi^2}{3}},\quad \nu(n +1)=\sqrt{\nu(n)^2+\frac{\pi^2}{3}}.
\end{equation}
  Observe that for $\nu(n)\geq 2$,
\begin{align*}
	&\nu(n-1)=\nu(n)-\frac{\pi^2}{6\nu(n)}-\frac{\pi^4}{72\nu(n)^3}
-\frac{\pi^6}{432\nu(n)^5}-\frac{5\pi^8}{10368\nu(n)^7}+
O\left(\frac{1}{\nu(n)^8}\right),\\[5pt]	&\nu(n+1)=\nu(n)+\frac{\pi^2}{6\nu(n)}-\frac{\pi^4}{72\nu(n)^3}+\frac{\pi^6}{432\nu^5(n)}
-\frac{5\pi^8}{10368\nu(n)^7}+O\left(\frac{1}{\nu(n)^8}\right),
\end{align*}
so it is readily checked that  for $\nu(n)\geq 3$,
\begin{align}
	d_v(n)<&\nu(n-1)<u_v(n), \label{xkn-1}\\[5pt]
	\bar{d}_v(n)<&\nu(n+1)<\bar{u}_v(n), \label{xkn+1}
\end{align}
where
\begin{align}
\begin{aligned}\label{wylabel}
	d_v(n)&=\nu(n)-\frac{\pi^2}{6\nu(n)}-\frac{\pi^4}{72\nu(n)^3}
	-\frac{\pi^6}{432\nu(n)^5}-\frac{5\pi^8}{5184\nu(n)^7},\\
	u_v(n)&=\nu(n)-\frac{\pi^2}{6\nu(n)}-\frac{\pi^4}{72\nu(n)^3}
	-\frac{\pi^6}{432\nu(n)^5},\\
\bar{d}_v(n)&=\nu(n)+\frac{\pi^2}{6\nu(n)}-\frac{\pi^4}{72\nu(n)^3}+\frac{\pi^6}{432\nu(n)^5}
-\frac{5\pi^8}{5184\nu(n)^7},\\
	\bar{u}_v(n)&=\nu(n)+\frac{\pi^2}{6\nu(n)}-\frac{\pi^4}{72\nu(n)^3}
	+\frac{\pi^6}{432\nu(n)^5}.
\end{aligned}
\end{align}

With  \eqref{xkn-1} and \eqref{xkn+1} in hands, we are now in a position to bound   $\exp(\nu(n-1)\break$
$+\nu(n+1)-2\nu(n))$, $L(n)$ and $R(n)$ in terms of $\nu(n)$.

We first estimate $\exp( \nu(n-1)+\nu(n+1)-2\nu(n))$. Applying \eqref{xkn-1} and \eqref{xkn+1}, we find that for $\nu(n) \geq 3$,
\[-\frac{ \pi^4}{36 \nu(n)^3}-\frac{5  \pi^8}{2592 \nu(n)^7}<\nu(n-1)+\nu(n+1)-2\nu(n)<-\frac{ \pi^4}{36 \nu(n)^3}.\]
It follows that for $\nu(n) \geq 3$,
\begin{equation}\label{eq-exp-1}
\exp( {  \nu(n-1)+\nu(n+1)-2\nu(n)})<\exp\left({-\frac{ \pi^4}{36 \nu(n)^3}}\right)
\end{equation}
and
\begin{equation}\label{eq-exp-2}
\exp( {  \nu(n-1)+\nu(n+1)-2\nu(n)})>\exp\left( {-\frac{ \pi^4}{36 \nu(n)^3}-\frac{5  \pi^8}{2592 \nu(n)^7}}\right).
\end{equation}
Note that for $s<0$,
\begin{equation*}
1+s<e^s<1+s+s^2.
\end{equation*}
Hence we get that
\begin{equation}\label{eq-exp-3}
\exp\left(-\frac{ \pi^4}{36 \nu(n)^3}\right)<1-\frac{ \pi^4}{36 \nu(n)^3}+\frac{ \pi^8}{1296 \nu(n)^6}
\end{equation}
and
\begin{equation}\label{eq-exp-4}
\exp\left(-\frac{ \pi^4}{36 \nu(n)^3}-\frac{5  \pi^8}{2592 \nu(n)^7}\right)>1-\frac{ \pi^4}{36 \nu(n)^3}-\frac{5  \pi^8}{2592 \nu(n)^7}.
\end{equation}
Combining \eqref{eq-exp-1} and \eqref{eq-exp-3}  yields that for $\nu(n) \geq 3$,
\begin{equation}\label{eq-exp-5}
 \exp\left({  \nu(n-1)+\nu(n+1)-2\nu(n)}\right)<1-\frac{ \pi^4}{36 \nu(n)^3}+\frac{ \pi^8}{1296 \nu(n)^6} .
\end{equation}
Using \eqref{eq-exp-2} together with  \eqref{eq-exp-4}, we get that for $\nu(n) \geq 3$,
\begin{equation}\label{eq-exp-6}
 \exp\left({  \nu(n-1)+\nu(n+1)-2\nu(n)}\right)>1-\frac{ \pi^4}{36 \nu(n)^3}-\frac{5  \pi^8}{2592 \nu(n)^7}.
\end{equation}

Next, we   estimate $L(n)$ and $R(n)$. Let
\begin{align}
		{P}_l(n)
		&=\frac{1}{\nu(n-1)^6\nu(n
+1)^{6}}\left(\nu(n-1)^6-\frac{3}{8}
		\nu(n-1)^4u_v(n)-\frac{15}{128} \nu(n-1)^4\right. \nonumber \\[6pt]
		&\phantom{=\;\;}\left.\qquad-\frac{105}{1024} \nu(n-1)^2u_v(n)-\frac{4725}{32768} \nu(n-1)^2-\frac{72765}{262144} u_v(n)-31\right)\nonumber \\[6pt]
		&\phantom{=\;\;}\quad\times\left(\nu(n+1)^6-\frac{3}{8}\nu(n+1)^4\bar{u}_v(n)-\frac{15}{128} \nu(n+1)^4-\frac{105}{1024} \nu(n+1)^2\bar{u}_v(n)\right.\nonumber \\[6pt]
		&\phantom{=\;\;}\left.\qquad-\frac{4725}{32768} \nu(n+1)^2-\frac{72765}{262144} \bar{u}_v(n)-31\right) \label{def-pl}
\end{align}
and
\begin{align}
	 {P}_r(n)
	&=\frac{1}{\nu(n-1)^6\nu(n
+1)^{6}}\left(\nu(n-1)^6-\frac{3}{8}
	\nu(n-1)^4d_v(n)-\frac{15}{128} \nu(n-1)^4\right. \nonumber\\[6pt]
	&\phantom{=\;\;}\left.\qquad-\frac{105}{1024} \nu(n-1)^2d_v(n)-\frac{4725}{32768} \nu(n-1)^2-\frac{72765}{262144} d_v(n)+31\right) \nonumber\\[6pt]
	&\phantom{=\;\;}\quad\times\left(\nu(n+1)^6-\frac{3}{8}\nu(n+1)^4\bar{d}_v(n)-\frac{15}{128} \nu(n+1)^4-\frac{105}{1024} \nu(n+1)^2\bar{d}_v(n)\right. \nonumber \\[6pt]
	&\phantom{=\;\;}\left.\qquad-\frac{4725}{32768} \nu(n+1)^2-\frac{72765}{262144} \bar{d}_v(n)+31\right). \label{def-pr}
\end{align}
Applying \eqref{xkn-1} and \eqref{xkn+1} into \eqref{defi-tl} and \eqref {defi-hl},   it is not difficult to check that for $\nu(n) \geq 3$,
\begin{equation}\label{eq-tl-1}
	L(n)\geq\frac{{P}_l(n)}{ \left(E_{I}(\nu(n))+\frac{31}{\nu(n)^6}\right)^2}
\end{equation}
and
\begin{equation}\label{eq-hl-1}
	R(n)\leq\frac{{P}_r(n)}{  \left(E_{I}(\nu(n))-\frac{31}{\nu(n)^6}\right)^2}.
\end{equation}
To bound $L(n)$ and $R(n)$ in terms of $\nu(n)$, we shall show that for $\nu(n) \geq 60$,
\begin{align}\label{eq-tl-2}
\frac{{P}_l(n)}{ \left(E_{I}(\nu(n))+\frac{31}{\nu(n)^6}\right)^2}&
\geq 1-\frac{ \pi^4}{32  \nu(n)^5}-\frac{129}{\nu(n)^6}\nonumber\\[6pt]
&=\frac{32\nu(n)^6-\pi^4\nu(n)-4128}{32\nu(n)^6}
\end{align}
and
\begin{align}\label{eq-hl-2}
	\frac{{P}_r(n)}{  \left(E_{I}(\nu(n))-\frac{31}{\nu(n)^6}\right)^2}
	& \leq 1-\frac{ \pi^4}{32  \nu(n)^5}+\frac{121}{  \nu(n)^6} \nonumber\\[6pt]
	&=\frac{32\nu(n)^6-\pi^4\nu(n)+3872}{32\nu(n)^6},
\end{align}
which are equivalent to showing that for $\nu(n) \geq 60$,
\begin{align}\label{eq-tl-3}
	32\nu(n)^6{P}_l(n)-\left(32\nu(n)^6-\pi^4\nu(n)-4128\right)
	\left(E_{I}(\nu(n))+\frac{31}{\nu(n)^6}\right)^2\geq 0
\end{align}
and
\begin{align}\label{eq-hl-3}
\left(32\nu(n)^6-\pi^4\nu(n)+3872\right)\left(E_{I}(\nu(n))-\frac{31}{\nu(n)^6}\right)^2-32\nu(n)^6{P}_r(n)\geq 0.
\end{align}
Substituting   \eqref{eq-x-Rel} and \eqref{wylabel} into \eqref{def-pl} and \eqref{def-pr}, we find that
\begin{align}\label{eq-tl-4}
	&32\nu(n)^6{P}_l(n)-\left(32\nu(n)^6-\pi^4\nu(n)-4128\right)
	\left(E_{I}(\nu(n))+\frac{31}{\nu(n)^6}\right)^2\nonumber \\[5pt]
	& =\frac{\sum_{j=0}^{26} a_j \nu(n)^j }{\nu(n)^{14}\nu(n-1)^6\nu(n+1)^6}
\end{align}
and
\begin{align}\label{eq-hl-4}
&\left(32\nu(n)^6-\pi^4\nu(n)+3872\right)\left(E_{I}(\nu(n))-\frac{31}{\nu(n)^6}\right)^2-32\nu(n)^6{P}_r(n)\nonumber\\[5pt]
&=\frac{\sum_{j=0}^{26} b_j \nu(n)^j}{\nu(n)^{14}\nu(n-1)^6\nu(n+1)^6},
\end{align}
where $a_j$ and $b_j$ are real numbers. Here we just list the values of $a_{24}$,  $a_{25}$,  $a_{26}$,  $b_{24}$,  $b_{25}$ and  $b_{26}$:
\[a_{24}=78-\frac{175 \pi ^4}{64},~~
	a_{25}=-1608 -\frac{19 \pi^4}{16},~~a_{26}=160-\frac{4 \pi ^4}{3},
\]
\[	b_{24}=102 + \frac{175 \pi ^4}{64},~~ b_{25}=\frac{19 \pi ^4}{16}-1416,~~
	b_{26}=\frac{4 \pi ^4}{3}-96.
\]
It can be readily checked that for any $0\leq j\leq 23$ and   $\nu(n)\geq 27$,
\[-|a_j|\nu(n)^j\geq -|a_{24}|\nu(n)^{24}\]
and
\[-|b_j|\nu(n)^j\geq -|b_{24}|\nu(n)^{24}.\]
It follows that for $\nu(n)\geq 27$,
\begin{align*}
	\sum_{j=0}^{26} a_j \nu(n)^j & \geq-\sum_{j=0}^{24}\left|a_j\right| \nu(n)^j+a_{25} \nu(n)^{25}+a_{26} \nu(n)^{26} \nonumber\\[6pt]
	& \geq-25\left|a_{24}\right| \nu(n)^{24}+a_{25} \nu(n)^{25}+a_{26} \nu(n)^{26}
\end{align*}
and
\begin{align*}
	\sum_{j=0}^{26} b_j \nu(n)^j & \geq-\sum_{j=0}^{24}\left|b_j\right| \nu(n)^j+b_{25} \nu(n)^{25}+b_{26} \nu(n)^{26} \nonumber\\[6pt]
	& \geq-25\left|b_{24}\right| \nu(n)^{24}+b_{25} \nu(n)^{25}+b_{26} \nu(n)^{26}.
\end{align*}
Moreover,  one can easily check that for $\nu(n) \geq 60$,
$$
-25\left|a_{24}\right| \nu(n)^{24}+a_{25} \nu(n)^{25}+a_{26} \nu(n)^{26} \geq 0
$$
and
$$
-25\left|b_{24}\right| \nu(n)^{24}+b_{25} \nu(n)^{25}+b_{26} \nu(n)^{26}\geq 0 .
$$
Hence \eqref{eq-tl-3} and \eqref{eq-hl-3} hold for $\nu(n) \geq 60$, and so  \eqref{eq-tl-2} and \eqref{eq-hl-2} hold for $\nu(n) \geq 60$. Substituting \eqref{eq-tl-2} into \eqref{eq-tl-1}, we obtain that for $\nu(n) \geq 60$,
\begin{align}\label{eq-tl-5}
	L(n)
	\geq 1-\frac{ \pi^4}{32  \nu(n)^5}-\frac{129}{\nu(n)^6}.
\end{align}
Plugging \eqref{eq-hl-2} into \eqref{eq-hl-1}, we get that for $\nu(n) \geq 60$,
\begin{align}\label{eq-hl-5}
	R(n)
	\leq 1-\frac{ \pi^4}{32  \nu(n)^5}+\frac{121}{\nu(n)^6}.
\end{align}
Applying \eqref{eq-exp-5}, \eqref{eq-exp-6}, \eqref{eq-tl-5} and \eqref{eq-hl-5} into \eqref{eq-I-1} and \eqref{eq-I-2}, we are lead to \eqref{eq-lem-B1} and \eqref{eq-lem-B2}. This completes the proof.
\qed\\

	\section{Proofs of Theorem \ref{eq-q-thm} and Theorem \ref{lem-M-1}}

To prove Theorem \ref{eq-q-thm}, we first derive  an asymptotic formula for $q(n)$  with an explicit bound  by specializing  an asymptotic formulas for $\eta$-quotients $G(q)$ due to  Chern \cite{Chern-2019}. Define
	\begin{eqnarray}
		G(q)=G(e^{2\pi i \tau}):=\prod_{r=1}^R(q^{m_r};q^{m_r})_{\infty}^{\delta_r},
		\label{eq-1-1}
	\end{eqnarray}
	where $\mathbf{m}=(m_1,\ldots,m_R)$ is a sequence of  $R$ distinct positive integers and $\mathbf{\delta}=(\delta_1,\ldots,\delta_R)$ is a sequence of $R$ non-zero integers.
	Here and throughout this paper,  we have adopted the  standard notation on  $q$-series \cite{Andrews-1998}.
	\[(a;q)_n=\prod_{j=0}^{n-1}(1-aq^j) \quad \text{and} \quad (a;q)_\infty=\prod_{j=0}^\infty (1-aq^j).\]
	
	In order to state Chern's result, we need a few preliminary definitions. Assume that $h$ and $j$ are positive integers with ${\rm gcd}(h, j) =1$,  set
\[\Delta_1=-\frac{1}{2}\sum_{r=1}^R\delta_r,\qquad \Delta_2=\sum_{r=1}^Rm_r\delta_r,\]
\[\Delta_3(k)=-\sum_{r=1}^R\frac{\delta_r\gcd^2(m_r,k)}{m_r},\qquad \Delta_4(k)=\prod_{r=1}^R\left(\frac{m_r}{\gcd(m_r,k)}\right)^{-\frac{\delta_r}{2}},\]

\begin{equation}\label{defi-A}
		\hat{A}_{k}(n)=\sum_{0\leq h< k\atop {\rm gcd}(h,k)=1}\exp\left(-\frac{2\pi nh i }{k}-\pi i\sum_{r=1}^R\delta_rs\left(\frac{m_rh}{\gcd(m_r,k)},\frac{k}{\gcd(m_r,k )}\right)\right),
	\end{equation}
	where  $s(h,j)$ is the  Dedekind sum defined by
	\[s(h,j)=\sum_{r=1}^{j-1}\left(\frac{r}{j}-\left[\frac{r}{j}\right]-\frac{1}{2}\right)\left(\frac{h r}{j}-\left[\frac{h r}{j}\right]-\frac{1}{2}\right).\]
	
Let $L=\mathrm {lcm} (m_1,\ldots,m_R)$. We divide the set  $\{1,2\cdots,L\}$ into two disjoint subsets:
\[\mathcal{L}_{>0}:=\{1\le l\le L:\Delta_3(l)>0\},\]
\[\mathcal{L}_{\le0}:=\{1\le l\le L:\Delta_3(l)\le0\}.\]
If we write
\[G(q)=\sum_{n\geq 0} g(n)q^n.\]
Chern \cite{Chern-2019} obtained an asymptotic formula for $g(n)$ with $\frac{1}{2}\sum_{r=1}^R\delta_r\geq 0$.

Define
\begin{align}
	\mathbb{E}_{\Delta_1}(s):=
	\left\{\begin{array}{ll}
		1,~~~~~~~~~~~~~~~~~~&\Delta_1=0,\\[5pt]
		2\sqrt{s},~~~~~~~~~~~&\Delta_1=-\frac{1}{2},\\[5pt]
		s\log(s+1),~~~~~~~&\Delta_1=-1,\\[5pt]
		s^{-2\Delta_1-1}\zeta(-\Delta_1),~~~~~~~&\mbox{otherwise},\\[5pt]
	\end{array}\right.
	\label{eq-E-N}
\end{align}
where $\zeta(\cdot)$ is Riemann zeta-function.

	\begin{thm}[Chern]\label{thm-Chern} 	If  $\Delta_1\leq0$ and the inequality
		\begin{equation}\label{eq-min}
			\min_{1\le r\le R}\left(\frac{\gcd^2(m_r,l)}{m_r}\right)\ge \frac{\Delta_3(l)}{24}
		\end{equation}
		holds for all  $1\le l \le L$, then for positive integers $N$ and  $n>-\frac{\Delta_2}{24}$, we have
	\begin{align}
		g(n)=E(n)+&\sum_{l\in \mathcal{L}_{>0}}2\pi \Delta_4(l)\left(\frac{24n+\Delta_2}{\Delta_3(l)}\right)^{-\frac{\Delta_1+1}{2}}\nonumber\\[5pt]
		&\quad\times\sum_{\substack{1\le k\le N\\k\equiv_L l }}\frac{I_{-\Delta_1-1}\left(\frac{\pi}{6k}\sqrt{\Delta_3(l)(24n+\Delta_2)}\right)}{k}\hat{A}_{k}(n),
	\end{align}
where
\begin{align}
	|E(n)|
	&\le\frac{2^{-\Delta_1}\pi^{-1}N^{-\Delta_1+2}}{n+\frac{\Delta_2}{24}}\exp\left(
 2\pi \left(n+\frac{\Delta_2}{24}\right)N^{-2}\right)\sum_{l\in \mathcal{L}_{>0}}\Delta_4(l)\exp\left(\frac{\Delta_3(l)\pi}{3}\right)\nonumber\\[5pt]
	&\quad+2\exp\left(2\pi\left(n+\frac{\Delta_2}{24}\right)N^{-2}\right)\mathbb{E}_{\Delta_1}(N)\nonumber\\[5pt]
	&\quad\times\left(\sum_{1\le l\le L}\Delta_4(l)\exp\left(\frac{\pi \Delta_3(l)}{24}+\sum_{r=1}^R\frac{|\Delta_r|\exp\left(-\pi\gcd^2(m_r,l)/m_r\right)}{\left(1-\exp\left(-\pi\gcd^2(m_r,l)/m_r\right)\right)^2}\right)\right.\nonumber\\[5pt]
       &\phantom{=\;\;}\left.\quad-\sum_{l\in \mathcal{L}_{>0}}\Delta_4(l)\exp\left(\frac{\pi \Delta_3(l)}{24}\right)\right)\nonumber,
\end{align}
and $I_{\nu}(s)$ is the $\nu$-th modified Bessel function of the first kind.	
	%

	\end{thm}

We are now in a position to prove Theorem \ref{eq-q-thm} by means of  Theorem \ref{thm-Chern}.

\noindent {\it Proof of Theorem \ref{eq-q-thm}.}   Recall that	
	\begin{align*}
		\sum_{n=0}^{\infty}q(n) q^{n} =  \frac{(q^2;q^2)_{\infty}}{(q;q)_{\infty}},
	\end{align*}
so we have
	$\mathbf{m}=(1,2)$ and $\mathbf{\delta}=(-1,1)$. It is straightforward to compute that $\Delta_1=0,$ and $\Delta_2=1$. We also have $L=2$. The values of $\Delta_3(l)$ and $\Delta_4(l)$ for $1\le l\le L$
are listed in Table \ref{tab:1}. Hence $\mathcal{L}_{>0}=\{1\}$.  It can be readily checked that \eqref{eq-min}  is always true for $1\leq l \leq 2$.
\begin{table}[htbp]
	\centering
	
	\caption{The values of $\Delta_3(l)$ and $\Delta_4(l)$ for $1\le l\le 2$.}
	
	\label{tab:1}       
	\begin{tabular}{ccc}
		
		\hline\noalign{\smallskip}
		
		$l$ & 1 & 2  \\
		
		\noalign{\smallskip}\hline\noalign{\smallskip}
		
		$\Delta_3(l)$  & $ \frac{1}{2}$ &  $-1$  \\
		
		$\Delta_4(l)$ & $\frac{\sqrt{2}}{2}$& $1$ \\
		
		\noalign{\smallskip}\hline
		
	\end{tabular}
\end{table}

Hence, by Theorem \ref{thm-Chern}, we have
\begin{equation*}
	q(n)=E(n)+\frac{   \sqrt{2}\pi^2}{12 \nu(n)}\sum_{\substack{1\le k\le N \\2\nmid k }}I_1\left(\frac{\nu(n)}{k}\right)\frac{\hat{A}_k(n)}{k},
\end{equation*}
where $\nu(n)$ is defined as in \eqref{defi-x}, and
\begin{align}
|E(n)|&\le\frac{\pi^{-1}N^{2}}{n+\frac{1}{24}}\exp\left(2\pi\left(n+\frac{1}{24}\right)N^{-2}\right)\cdot\frac{\sqrt{2}}{2}\exp\left(\frac{\pi}{6}\right)+2\exp\left(2\pi\left(n+\frac{1}{24}\right)N^{-2}\right)\nonumber\\[5pt]
&\quad\times\left\{ \frac{\sqrt{2}}{2}\exp\left(\frac{\pi }{48}+\frac{\exp(-\pi)}{(1-\exp(-\pi))^2}+\frac{\exp(-\pi/2)}{(1-\exp(-\pi/2))^2}\right)
\right.\nonumber\\[5pt]
&\phantom{=\;\;}\left.
\quad+\exp\left(-\frac{\pi }{24}+\frac{\exp(-\pi)}{(1-\exp(-\pi))^2}+\frac{\exp(-2\pi)}{(1-\exp(-2\pi))^2}\right)-\frac{\sqrt{2}}{2}\exp\left(\frac{\pi }{48}\right)\right\}\nonumber.
\end{align}	
Assume that $N=\left\lfloor \nu(n)\right\rfloor$,  then
\begin{align}
	|E(n)|&\le\frac{\sqrt{2}\pi }{6}\frac{\left\lfloor \nu(n)\right\rfloor^2}{\nu(n)^2}\exp\left(\frac{6\nu(n)^2}{\pi\left\lfloor \nu(n)\right\rfloor^{2}}+\frac{\pi}{6}\right)\nonumber\\[5pt]
	&\quad+2\exp\left(\frac{6\nu(n)^2}{\pi\left\lfloor \nu(n)\right\rfloor^{2}}\right)\left\{ \frac{\sqrt{2}}{2}\exp\left(\frac{\pi }{48}+\frac{\exp(-\pi)}{(1-\exp(-\pi))^2}+\frac{\exp(-\pi/2)}{(1-\exp(-\pi/2))^2}\right)
	\right.\nonumber\\[5pt]
	&\phantom{=\;\;}\left.
	\quad+\exp\left(-\frac{\pi }{24}+\frac{\exp(-\pi)}{(1-\exp(-\pi))^2}+\frac{\exp(-2\pi)}{(1-\exp(-2\pi))^2}\right)-\frac{\sqrt{2}}{2}\exp\left(\frac{\pi }{48}\right)\right\}\nonumber.
\end{align}
Using the following two inequalities:
\[\frac{\left\lfloor \nu(n)\right\rfloor^2}{\nu(n)^2}\le 1  \quad \text{and} \quad
\frac{\nu(n)^2}{\left\lfloor \nu(n)\right\rfloor^2}<\frac{\nu(n)^2}{(\nu(n)-1)^2}< 2\, \quad\text{for }   \nu(n)\ge 4,\]
we deduce  that
\begin{align*}
	|E(n)|&\le\frac{\sqrt{2}\pi }{6}\exp\left(\frac{12}{\pi}+\frac{\pi}{6}\right)\nonumber\\[5pt]
 &\quad+2\exp\left(\frac{12}{\pi}\right)\left\{ \frac{\sqrt{2}}{2}\exp\left(\frac{\pi }{48}+\frac{\exp(-\pi)}{(1-\exp(-\pi))^2}+\frac{\exp(-\pi/2)}{(1-\exp(-\pi/2))^2}\right)
	\right.\nonumber\\[5pt]
	&\phantom{=\;\;}\left.
	\quad+\exp\left(-\frac{\pi }{24}+\frac{\exp(-\pi)}{(1-\exp(-\pi))^2}+\frac{\exp(-2\pi)}{(1-\exp(-2\pi))^2}\right)-\frac{\sqrt{2}}{2}\exp\left(\frac{\pi }{48}\right)\right\}\le 173.
\end{align*}	
Thus, we conclude that for $\nu(n)\geq 4$,
 \begin{equation}\label{formula-q(n)}
	q(n)=E(n)+\frac{   \sqrt{2}\pi^2}{12 \nu(n)}\sum_{\substack{1\le k\le \left\lfloor \nu(n)\right\rfloor \\2\nmid k }}I_1\left(\frac{\nu(n)}{k}\right)\frac{\hat{A}_k(n)}{k},
\end{equation}
where   $|E(n)|\le 173$.

Observing that $\hat{A}_1(n)=1$, so by \eqref{formula-q(n)}, we have
\[q(n)=\frac{   \sqrt{2}\pi^2}{12 \nu(n)}I_{1}(\nu(n))+R(n),\]
where
\begin{equation}\label{eq-Rq}
R(n)=E(n)+\frac{   \sqrt{2}\pi^2}{12 \nu(n)}\sum_{\substack{3\le k\le \left\lfloor \nu(n)\right\rfloor \\2\nmid k }}I_1\left(\frac{\nu(n)}{k}\right)\frac{\hat{A}_k(n)}{k}.
\end{equation}
Hence, in order to show Theorem \ref{eq-q-thm}, it suffices to show that for $\nu(n)\geq 21$,
	\begin{equation}\label{bound-R}
		|R(n)|\leq \frac{\sqrt{3} \pi^{\frac{3}{2}}}{  6\nu(n)^{\frac{1}{2}}}  \exp\left(\frac{\nu(n)}{3}\right).
	\end{equation}
By the definition of  $\hat{A}_k(n)$, we derive that  for any $n\geq 0$ and $k\geq 1$,
\begin{equation*}\label{eq-A-ineq}
	|\hat{A}_k(n)|\leq k,
\end{equation*}
since $|e^{2\pi s i}|=1$ for any $s\in \mathbb{R}$.
It yields that
\begin{align*}
	\left|\frac{   \sqrt{2}\pi^2}{12 \nu(n)}\sum_{\substack{3\le k\le \left\lfloor \nu(n)\right\rfloor \\2\nmid k }}I_1\left(\frac{\nu(n)}{k}\right)\frac{\hat{A}_k(n)}{k}\right|&\leq\frac{ \sqrt{2}\pi^2}{12\nu(n)}\sum_{3\le k\le \left\lfloor \nu(n)\right\rfloor\atop 2\nmid k } I_{1}\left(\frac{\nu(n)}{k}\right)\\[6pt]
	&\leq\frac{ \sqrt{2}\pi^2}{12\nu(n)}\frac{\left\lfloor \nu(n)\right\rfloor}{2}I_1\left(\frac{\nu(n)}{3}\right)\\[6pt]
	&\leq \frac{ \sqrt{2}\pi^2}{24}I_1\left(\frac{\nu(n)}{3}\right).
\end{align*}
Invoking \eqref{ineq-BKRT}, we derive that
\begin{align}\label{eq-RI}
	\left|\frac{   \sqrt{2}\pi^2}{12 \nu(n)}\sum_{\substack{3\le k\le N \\2\nmid k }}I_1\left(\frac{\nu(n)}{k}\right)\frac{\hat{A}_k(n)}{k}\right|&\leq   \frac{\sqrt{3} \pi^{\frac{3}{2}}}{ 12\nu(n)^{\frac{1}{2}}}  \exp\left(\frac{\nu(n)}{3}\right).
\end{align}
Substituting  \eqref{eq-RI} into \eqref{eq-Rq}, we get that for $\nu(n)\geq 4$,
\begin{equation}\label{eq-R}
	|R(n)|\leq 173+\frac{\sqrt{3} \pi^{\frac{3}{2}}}{  12\nu(n)^{\frac{1}{2}}}  \exp\left(\frac{\nu(n)}{3}\right) .
\end{equation}
We proceed to show that for $\nu(n)\geq21 $,
\begin{equation}\label{eq-r}
	173<\frac{\sqrt{3} \pi^{\frac{3}{2}}}{  12\nu(n)^{\frac{1}{2}}}  \exp\left(\frac{\nu(n)}{3}\right).
\end{equation}
Define
\begin{equation*}
r(s):=	\frac{692\sqrt{3}}{ \pi^{\frac{3}{2}}  }s^{\frac{1}{2}} \exp\left(-\frac{s}{3}\right).
\end{equation*}
It is evident that
\begin{align*}
	r'(s)=	\frac{346 }{\sqrt{3} \pi ^{\frac{3}{2}}s^{\frac{1}{2}}}(-2 s+3)\exp\left(-\frac{s}{3}\right).
\end{align*}
Since  $r'(s)\leq 0$ when $s\geq\frac{3}{2}$, we deduce that $r(s)$ is decreasing when $s \geq\frac{3}{2}$. This implies that
\[r(\nu(n))\leq r(21)<1\]
for $\nu(n)\geq 21$. So the inequality \eqref{eq-r} is valid.	
Applying \eqref{eq-r} to \eqref{eq-R}, we are led to \eqref{bound-R}. This completes the proof. \qed

We conclude this section with the proof of Theorem \ref{lem-M-1} by employing Theorem  \ref{eq-q-thm}.

{\noindent \it Proof of Theorem \ref{lem-M-1}.}
Define
\begin{equation}\label{defi-Gk}
	G(n):=\frac{\frac{\sqrt{3} \pi^{\frac{3}{2}}}{ 6\nu(n)^{\frac{1}{2}}}  \exp\left(\frac{\nu(n)}{3}\right)} {\frac{\sqrt{2}\pi^2}{12\nu(n)}I_1(\nu(n))}= \sqrt{\frac{6\nu(n)}{\pi}}\cdot\frac{\exp\left(\frac{\nu(n)}{3}\right)}{I_1(\nu(n))}.
\end{equation}
Thanks to  Theorem \ref{eq-q-thm}, we have
\begin{equation*}
	M(n)(1-G(n))\leq q(n)\leq M(n)(1+G(n)).
\end{equation*}
To show \eqref{eq-main}, it is enough to prove that for $\nu(n)\geq38$,
\begin{equation}
	G(n)\leq \frac{1}{\nu(n)^6}.
\end{equation}
Using Lemma \ref{thm-bessel-1}, we find that for $s\geq 26$,
\begin{equation*}
	I_1(s)\geq \frac{e^s}{\sqrt{2\pi s}}\left(1-\frac{3}{8s}-\frac{15}{128s^2}-\frac{105}{1024s^3}
-\frac{4725}{32768s^4}-\frac{72765}{262144s^5}-\frac{31}{s^6}\right).
\end{equation*}
Note that for $s\geq 4$
\begin{equation*}
	\frac{1}{8 s}-\frac{15}{128s^2}-\frac{105}{1024 s^3}-\frac{4725}{32768s^4}-\frac{72765}{262144s^5}-\frac{31}{s^6}\geq 0,
\end{equation*}
so for $s\geq26$,
\begin{equation}\label{emitate-bessel}
	I_1(s) \geq \frac{e^{s}}{\sqrt{2\pi s}} \left(1-\frac{1}{2 s}\right).
\end{equation}
Substituting \eqref{emitate-bessel} into \eqref{defi-Gk},  we derive that
 for $\nu(n)\geq 26$,
\begin{align}\label{eq-S}
	G(n)\leq\frac{2\sqrt{3}\nu(n)}{1-\frac{1}{2\nu(n)}}\exp\left(-\frac{2\nu(n)}{3}\right).	
\end{align}
Based on the following observation:
\begin{align*}
	\left(1-\frac{1}{2 \nu(n)}\right)\left(1+\frac{1}{ \nu(n)}\right)&=1+\frac{1}{2\nu(n)^2}\left(\nu(n)-1\right)\geq 1 \quad  \text{for} \quad \nu(n)\geq 1,
\end{align*}
  we find that  \eqref{eq-S} can be further bounded by
\begin{equation}\label{eq-G(x)}
	G(n)    \leq 2\sqrt{3}\nu(n)\left(1+\frac{1}{ \nu(n)}\right)\exp\left(-\frac{2\nu(n)}{3}\right)	.
\end{equation}
We claim that  for $\nu(n)\geq 43$,
\begin{equation}\label{eq-f(x)}
2\sqrt{3}\exp\left(-\frac{2\nu(n)}{3}\right) \leq\frac{1}{2\nu(n)^7},
\end{equation}
which can be recast as
\[4\sqrt{3}\nu(n)^7\exp\left(-\frac{2\nu(n)}{3}\right) \leq 1.
\]
Define
\[L(s):=4\sqrt{3}s^7\exp\left(-\frac{2s}{3}\right).\]
Since
\begin{align*}
L'(s)=	4\sqrt{3} \exp\left(-\frac{2s}{3}\right) s^6 \left(- \frac{2}{3}s+7\right) \leq 0 \quad \text{for } s\geq\frac{21}{2},
\end{align*}
we find that  $L(s)$ is decreasing when $s \geq\frac{21}{2}$. It follows that   for $\nu(n)\geq43$,
\[L(\nu(n))=4\sqrt{3}\nu(n)^7\exp\left(-\frac{2\nu(n)}{3}\right)\leq L(43)<1,\]
and so \eqref{eq-f(x)} holds when $\nu(n)\geq 43$. Hence the claim is verified.

Applying \eqref{eq-f(x)} to \eqref{eq-G(x)}, we are led to
\begin{align*}
	G(n)\leq \nu(n)\left(1+\frac{1}{\nu(n)} \right) \cdot \frac{1}{2\nu(n)^7}<\frac{1}{\nu(n)^6}
\end{align*}
 for  $\nu(n)\geq43 $. This completes the proof. \qed

\section{Proof of Theorem \ref{thm-B}}

In this section, we give a proof of Theorem \ref{thm-B} with the aid of  Theorem \ref{lem-M-1} and  Lemma \ref{lem-A-1}.

\noindent{\it Proof of Theorem \ref{thm-B}.} Recall that
\begin{equation*}
	Q(n)=\frac{q(n-1)q(n+1)}{q(n)^2}.
\end{equation*}
Define
\begin{equation}\label{defi-An}
	A(n)=\frac{M(n-1)M(n+1)}{M(n)^2},
\end{equation}
where $M(n)$ is defined as in \eqref{defi-M}.
From Theorem \ref{lem-M-1}, we see that for $\nu(n)\geq 43$,
\begin{equation}\label{eq-B-1}
A(n)L_Q(n)\leq Q(n)\leq A(n)R_Q(n),
\end{equation}
where
\begin{equation}\label{defi-tC}
	L_Q(n)=\frac{\left(1-\frac{1}{\nu(n-1)^6}\right)\left(1-\frac{1}{\nu(n+1)^6}\right)}{\left(1+\frac{1}{\nu(n)^6}\right)^2}
\end{equation}
and
\begin{equation}\label{defi-hC}
	R_Q(n)=\frac{\left(1+\frac{1}{\nu(n-1)^6}\right)\left(1+\frac{1}{\nu(n+1)^6}\right)}{\left(1-\frac{1}{\nu(n)^6}\right)^2}.
\end{equation}
To obtain \eqref{thm-B-low}, we proceed to estimate $A(n)$,  $L_Q(n)$ and $R_Q(n)$ in terms of $\nu(n)$. We first consider $A(n)$.
 Substituting \eqref{defi-M} into \eqref{defi-An}, we find that
\begin{equation}\label{defi-An-gene}
A(n)=\frac{\nu(n)^2I_1(\nu(n-1))I_1(\nu(n+1))}{\nu(n-1)\nu(n+1)I_1(\nu(n))^2}.
\end{equation}
Applying Lemma \ref{lem-A-1} into \eqref{defi-An-gene}, we obtain that for $\nu(n)\geq 60$,
\begin{align}\label{eq-A-3}
	A(n)
	&\geq \frac{\nu(n)^3}{\sqrt{\nu(n-1)^3\nu(n+1)^3}}\left(1-\frac{ \pi^4}{36 \nu(n)^3}-\frac{5  \pi^8}{2592 \nu(n)^7}\right)\nonumber\\[6pt]
	&\quad\times\left(1-\frac{ \pi^4}{32  \nu(n)^5}-\frac{129}{\nu(n)^6}\right)
\end{align}
and
\begin{align}\label{eq-A-4}
A(n)
	&\leq \frac{\nu(n)^3}{\sqrt{\nu(n-1)^3\nu(n+1)^3}}\left(1-\frac{ \pi^4}{36 \nu(n)^3}+\frac{ \pi^8}{1296 \nu(n)^6}\right) \nonumber\\[6pt]
	&\quad\times\left(1-\frac{ \pi^4}{32  \nu(n)^5}+\frac{121}{\nu(n)^6}\right).
\end{align}
We claim that for $\nu(n)\geq 8$,
\begin{equation}\label{eq-A-5}
1+\frac{\pi ^4}{12 \nu(n)^4}+\frac{7 \pi ^8}{864 \nu(n)^8}\leq \frac{\nu(n)^3}{\sqrt{\nu(n-1)^3\nu(n+1)^3}}\leq 1+\frac{\pi ^4}{12 \nu(n)^4}+\frac{ \pi ^8}{123 \nu(n)^8},
\end{equation}
which is equivalent to
\begin{align}\left\{
	\begin{aligned}\label{eq-A-6}
		&\nu(n)^{12}-\nu(n-1)^6\nu(n+1)^6\left(1+\frac{\pi ^4}{12 \nu(n)^4}+\frac{7 \pi ^8}{864 \nu(n)^8}\right)^4\geq 0,\\[6pt]
		&\nu(n)^{12}-\nu(n-1)^6\nu(n+1)^6\left(1+\frac{\pi ^4}{12 \nu(n)^4}+\frac{\pi ^8}{123\nu(n)^8}\right)^4\leq 0.\\
	\end{aligned}\right.
\end{align}
Recall that
\begin{equation*}
	\nu(n-1)=\sqrt{\nu(n)^2-\frac{\pi^2}{3}} \quad \text{and} \quad  \nu(n +1)=\sqrt{\nu(n)^2+\frac{\pi^2}{3}}.
\end{equation*}
It can be calculated that
\begin{align}\label{eq-A-7}
	&\nu(n)^{12}-\nu(n-1)^6\nu(n+1)^6\left(1+\frac{\pi ^4}{12 \nu(n)^4}+\frac{7 \pi ^8}{864 \nu(n)^8}\right)^4\nonumber\\[6pt]
	&=\frac{\pi^{12}}{406239826673664 \nu(n)^{32}}\left(1340897918976 \nu(n)^{32}+27935373312 \pi ^4 \nu(n)^{28}\right.\nonumber\\[6pt]
	&\phantom{=\;\;}\left.\quad+1551965184 \pi ^8 \nu(n)^{24}-1551965184 \pi ^{12} \nu(n)^{20}-60816096 \pi ^{16} \nu(n)^{16}\right.\nonumber\\[6pt]
	&\phantom{=\;\;}\left.\quad-3873177 \pi ^{20} \nu(n)^{12}+625779 \pi ^{24} \nu(n)^8+33957 \pi ^{28} \nu(n)^4+2401 \pi ^{32}\right)
\end{align}
and
\begin{align}\label{eq-A-8}
	&\nu(n)^{12}-\nu(n-1)^6\nu(n+1)^6\left(1+\frac{\pi ^4}{12 \nu(n)^4}+\frac{ \pi ^8}{123 \nu(n)^8}\right)^4\nonumber\\[6pt]
	&=-\frac{\pi^{8}}{42715740489984 \nu(n)^{32}}\left(4823367264 \nu(n)^{36}-141396118128 \pi ^4 \nu(n)^{32}\right.\nonumber\\[6pt]
	&\phantom{=\;\;}\left.	\quad-2942756919 \pi ^8 \nu(n)^{28}-175420755 \pi ^{12} \nu(n)^{24}+163918779 \pi ^{16} \nu(n)^{20}\right.\nonumber\\[6pt]
	&\phantom{=\;\;}\left.\quad+6413999 \pi ^{20} \nu(n)^{16}+418192 \pi ^{24} \nu(n)^{12}-66144 \pi ^{28} \nu(n)^8\right.\nonumber\\[6pt]
	&\phantom{=\;\;}\left.
	\quad-3584 \pi ^{32} \nu(n)^4-256 \pi ^{36}\right).
\end{align}
Note that for $\nu(n)\geq 4$,
\begin{align}\label{eq-A-9}
	&1551965184 \pi ^8 \nu(n)^{24}-1551965184 \pi ^{12} \nu(n)^{20}
\nonumber\\[6pt]
&\quad-60816096 \pi ^{16} \nu(n)^{16}-3873177 \pi ^{20} \nu(n)^{12}\geq 0,
\end{align}
and for $\nu(n)\geq 8$,
\begin{align}\left\{
	\begin{aligned}\label{eq-A-10}
		&4823367264 \nu(n)^{36}-141396118128 \pi ^4 \nu(n)^{32}\\[4pt]
		&\quad-2942756919 \pi ^8 \nu(n)^{28}-175420755 \pi ^{12} \nu(n)^{24}\geq 0,\\[6pt]
		&418192 \pi ^{24} \nu(n)^{12}-66144 \pi ^{28} \nu(n)^8-3584 \pi ^{32} \nu(n)^4-256 \pi ^{36}\geq 0.\\
	\end{aligned}\right.
\end{align}
Applying \eqref{eq-A-9} to \eqref{eq-A-7} and applying \eqref{eq-A-10} to \eqref{eq-A-8}, we find that \eqref{eq-A-6} holds for $\nu(n)\geq 8$, which implies \eqref{eq-A-5} holds for $\nu(n)\geq 8$, and so the claim is verified.  Substituting \eqref{eq-A-5} into \eqref{eq-A-3} and \eqref{eq-A-4}, we get that    for $\nu(n)\geq60$,
	
	\begin{align}\label{eq-A-1}
	A(n)
		&\geq \left(1+\frac{\pi ^4}{12 \nu(n)^4}+\frac{7 \pi ^8}{864 \nu(n)^8}\right)\left(1-\frac{ \pi^4}{36 \nu(n)^3}-\frac{5  \pi^8}{2592 \nu(n)^7}\right)\nonumber\\[6pt]
		&\quad\times\left(1-\frac{ \pi^4}{32  \nu(n)^5}-\frac{129}{\nu(n)^6}\right)
	\end{align}
	and
	\begin{align}\label{eq-A-2}
		A(n)
		&\leq \left(1+\frac{\pi ^4}{12 \nu(n)^4}+\frac{\pi ^8}{123 \nu(n)^8}\right)\left(1-\frac{ \pi^4}{36 \nu(n)^3}+\frac{ \pi^8}{1296 \nu(n)^6}\right) \nonumber\\[6pt]
		&\quad\times\left(1-\frac{ \pi^4}{32  \nu(n)^5}+\frac{121}{\nu(n)^6}\right).
	\end{align}
 We proceed to show that for $\nu(n)\geq 4$,
\begin{align}\label{eq-B-4}
	L_Q(n)\geq1-\frac{5}{\nu(n)^6} \quad \text{and} \quad
	R_Q(n)\leq1+\frac{5}{\nu(n)^6}.
\end{align}
Applying \eqref{eq-x-Rel} into \eqref{defi-tC} and \eqref{defi-hC},
  we find that
\begin{align*}
	L_Q(n)=\frac{\nu(n)^{12}\left(\left(\nu(n)^2+\frac{ \pi^2}{3}\right)^3-1\right)\left(\left(\nu(n)^2-\frac{ \pi^2}{3}\right)^3-1\right)}{\left(\nu(n)^6+1\right)^2\left(\nu(n)^4-\frac{ \pi^4}{9}\right)^3}
\end{align*}
and
\begin{align*}
	R_Q(n)=\frac{\nu(n)^{12}\left(\left(\nu(n)^2+\frac{ \pi^2}{3}\right)^3+1\right)\left(\left(\nu(n)^2-\frac{ \pi^2}{3}\right)^3+1\right)}{\left(\nu(n)^6-1\right)^2\left(\nu(n)^4-\frac{ \pi^4}{9}\right)^3}.
\end{align*}
Assume that
\begin{align*}
	\phi(s)=& 729 s^{24}-1215 \pi ^4 s^{20}+7290 s^{18}+81 \pi ^8 s^{16}-2187 \pi ^4 s^{14}+\left(3645-3 \pi ^{12}\right) s^{12}\\[6pt]
	&\quad
	+243 \pi ^8 s^{10}-1215 \pi ^4 s^8-9 \pi ^{12} s^6+135 \pi ^8 s^4-5 \pi ^{12}
\end{align*}
and
\begin{align*}
	\psi(s)=& 729 s^{24}-1215 \pi ^4 s^{20}-7290 s^{18}+81 \pi ^8 s^{16}+2187 \pi ^4 s^{14}+\left(3645-3 \pi ^{12}\right) s^{12}\\[6pt]
	&\quad-243 \pi ^8 s^{10}-1215 \pi ^4 s^8+9 \pi ^{12} s^6+135 \pi ^8 s^4-5 \pi ^{12}.
\end{align*}
It is not difficult to prove that
\begin{align}\label{eq-B-6}
L_Q(n)-\left(1-\frac{5}{\nu(n)^6}\right)=\frac{\phi\left(\nu(n)\right)}{\nu(n)^6\left(9 \nu(n)^4-\pi^4\right)^3\left(\nu(n)^6+1\right)^2}
\end{align}
and
\begin{align}\label{eq-B-7}
	R_Q(n)-\left(1+\frac{5}{\nu(n)^6}\right)=\frac{-\psi\left(\nu(n)\right)}{\nu(n)^6\left(9 \nu(n)^4-\pi^4\right)^3\left(\nu(n)^6-1\right)^2}.
\end{align}
Moreover, it is not difficult to show that
$\psi(s) \geq 0 $ for $s \geq 4$ and
\[
\phi(s) - \psi(s)=14580 s^{18}-4374 \pi ^4 s^{14}+486 \pi ^8 s^{10}-18 \pi ^{12} s^6>0 \]
for $s \geq 2$. Hence we derive that for $\nu(n) \geq 4$,
\begin{equation}\label{eq-B-8}
\phi\left(\nu(n)\right) >\psi\left(\nu(n)\right) \geq 0 .
\end{equation}
It follows that \eqref{eq-B-4} is valid.

Substituting \eqref{eq-A-1}, \eqref{eq-A-2} and   \eqref{eq-B-4} into \eqref{eq-B-1}, we derive that for $\nu(n)\ge 60$,
	\begin{align} \label{thm-B-pf-L}
	Q(n)
	&\geq \left(1+\frac{\pi ^4}{12 \nu(n)^4}+\frac{7 \pi ^8}{864 \nu(n)^8}\right)\left(1-\frac{ \pi^4}{36 \nu(n)^3}-\frac{5  \pi^8}{2592\nu(n)^7}\right) \nonumber \\[6pt]
	&\quad\times\left(1-\frac{ \pi^4}{32  \nu(n)^5}-\frac{129}{\nu(n)^6}\right)\left(1-\frac{5}{\nu(n)^6}\right)
\end{align}
and
\begin{align}\label{thm-B-pf-R}
	Q(n)
	&\leq \left(1+\frac{\pi ^4}{12 \nu(n)^4}+\frac{\pi ^8}{123 \nu(n)^8}\right)\left(1-\frac{ \pi^4}{36 \nu(n)^3}+\frac{ \pi^8}{1296 \nu(n)^6}\right) \nonumber \\[6pt]
	&\quad\times\left(1-\frac{ \pi^4}{32  \nu(n)^5}+\frac{121}{\nu(n)^6}\right)\left(1+\frac{5}{\nu(n)^6}\right).
\end{align}
To prove Theorem \ref{thm-B}, it is enough to show that  for  $\nu(n)\geq 67$,
\begin{align}\label{eq-B-9}
	&\left(1+\frac{\pi ^4}{12 \nu(n)^4}+\frac{7 \pi ^8}{864 \nu(n)^8}\right)\left(1-\frac{ \pi^4}{36 \nu(n)^3}-\frac{5  \pi^8}{2592 \nu(n)^7}\right)\nonumber\\[6pt]
	&\quad\times\left(1-\frac{ \pi^4}{32  \nu(n)^5}-\frac{129}{\nu(n)^6}\right)\left(1-\frac{5}{\nu(n)^6}\right)\nonumber\\[6pt]
	&> 1-\frac{\pi ^4}{36 \nu(n)^3}+\frac{\pi ^4}{12 \nu(n)^4}-\frac{\pi ^4}{32 \nu(n)^5}-\frac{135 }{\nu(n)^6}
\end{align}	
and	
\begin{align}\label{eq-B-10}
	&  \left(1+\frac{\pi ^4}{12 \nu(n)^4}+\frac{\pi ^8}{123 \nu(n)^8}\right)\left(1-\frac{ \pi^4}{36 \nu(n)^3}+\frac{ \pi^8}{1296 \nu(n)^6}\right) \nonumber\\[6pt]
	&\quad\times\left(1-\frac{ \pi^4}{32  \nu(n)^5}+\frac{121}{\nu(n)^6}\right)\left(1+\frac{5}{\nu(n)^6}\right)\nonumber\\[6pt]
	&< 1-\frac{\pi ^4}{36 \nu(n)^3}+\frac{\pi ^4}{12 \nu(n)^4}-\frac{\pi ^4}{32 \nu(n)^5}+\frac{126+\frac{\pi ^8}{1296}}{\nu(n)^6}.
\end{align}

We first show \eqref{eq-B-9}. Observe that
\begin{align}\label{eq-B-11}
	&\left(1+\frac{\pi ^4}{12 \nu(n)^4}+\frac{7 \pi ^8}{864 \nu(n)^8}\right)\left(1-\frac{ \pi^4}{36 \nu(n)^3}-\frac{5  \pi^8}{2592 \nu(n)^7}\right)\nonumber\\[6pt]
	&\quad\times\left(1-\frac{ \pi^4}{32  \nu(n)^5}-\frac{129}{\nu(n)^6}\right)\left(1-\frac{5}{\nu(n)^6}\right)\nonumber\\[6pt]
	&\quad-\left(1-\frac{\pi ^4}{36 \nu(n)^3}+\frac{\pi ^4}{12 \nu(n)^4}-\frac{\pi ^4}{32 \nu(n)^5}-\frac{135 }{\nu(n)^6}\right)\nonumber\\[6pt]
	&= \frac{1}{71663616 \nu(n)^{27}}\sum_{j=0}^{21}c_j\nu(n)^j,
\end{align}	
where $c_j$ are real numbers. Here we just list the values of $c_{19}$, $c_{20}$, $c_{21}$:
\begin{align}
	c_{19}=642816  \pi ^8,\quad 	c_{20}=-304128 \pi ^8,\quad	c_{21}=71663616 .\nonumber
\end{align}
  Clearly,
\begin{equation*}
	\sum_{j=0}^{21}c_j\nu(n)^j\geq -\sum_{j=0}^{19}|c_j|\nu(n)^j+c_{20}\nu(n)^{20}+c_{21}\nu(n)^{21}.
\end{equation*}
Moreover, it can be verified that for  $0\le j\le 18$ and  $\nu(n)\geq 4$,
\begin{equation*}
	-|c_j|\nu(n)^j\geq -|c_{19}|\nu(n)^{19}.
\end{equation*}
On the other hand, it's not difficult to check that for $\nu(n)\geq 67$,
\[c_{21}\nu(n)^2+c_{20}\nu(n)-20|c_{19}|> 0.\]
Assembling all these results above, we conclude that for $\nu(n)\ge 67$,
\begin{equation*}
	\sum_{j=0}^{21}c_j\nu(n)^j\geq \left(c_{21}\nu(n)^2+c_{20}\nu(n)-20|c_{19}|\right)\nu(n)^{19}> 0.
\end{equation*}	
This proves \eqref{eq-B-9}.

Similarly, to justify \eqref{eq-B-10},  we first note that
\begin{align}\label{eq-B-12}
	&  \left(1+\frac{\pi ^4}{12 \nu(n)^4}+\frac{\pi ^8}{123 \nu(n)^8}\right)\left(1-\frac{ \pi^4}{36 \nu(n)^3}+\frac{ \pi^8}{1296 \nu(n)^6}\right) \nonumber\\[6pt]
	&\quad\times\left(1-\frac{ \pi^4}{32  \nu(n)^5}+\frac{121}{\nu(n)^6}\right)\left(1+\frac{5}{\nu(n)^6}\right)\nonumber\\[6pt]
	&\quad-\left(1-\frac{\pi ^4}{36 \nu(n)^3}+\frac{\pi ^4}{12 \nu(n)^4}-\frac{\pi ^4}{32 \nu(n)^5}+\frac{126+\frac{\pi ^8}{1296}}{\nu(n)^6}\right)\nonumber\\[6pt]
	&=-\frac{1}{20404224\nu(n)^{26}}\sum_{j=0}^{19}d_j\nu(n)^j,
\end{align}
 where $d_j$ are real numbers. Here we also list the values of the last three coefficients:
\begin{align}
	d_{17}=53136 \pi ^8,\quad 	d_{18}=-183600 \pi ^8, \quad   d_{19}=47232 \pi ^8.
\end{align}
It's transparent that
\begin{equation}
\sum_{j=0}^{19}d_j\nu(n)^j\geq -\sum_{j=0}^{17}|d_j|\nu(n)^j+d_{18}\nu(n)^{18}+d_{19}\nu(n)^{19}.
\end{equation}
Moreover, it can be proved that for  $0\le j\le 16$ and $\nu(n)\geq2$,
\begin{equation*}
	-|d_j|\nu(n)^j\geq -|d_{17}|\nu(n)^{17}
\end{equation*}
and  for $\nu(n)\ge7$,
\[d_{19}\nu(n)^2+d_{18}\nu(n)-18|d_{17}|> 0.\]
Thus we conclude that for $\nu(n)\ge 67$,
\begin{equation*}
	\sum_{j=0}^{19}d_j\nu(n)^j\geq \left(d_{19}\nu(n)^2+d_{18}\nu(n)-18|d_{17}|\right)\nu(n)^{17}> 0,
\end{equation*}	
and so \eqref{eq-B-10} is valid.

Substituting \eqref{eq-B-9} and \eqref{eq-B-10} into \eqref{thm-B-pf-L} and  \eqref{thm-B-pf-R}, we arrive at   \eqref{thm-B-low}. This completes the proof.\qed

\section{Proof of Conjecture \ref{cp-1} }

In this section, we confirm Conjecture \ref{cp-1} with the aid of Theorem \ref{thm-B}. Before doing this, we first recall  the following lemma given by Jia \cite{Jia-2022}, which is useful  in  proof of  Conjecture \ref{cp-1} .

\begin{lem}[Jia]\label{lem-Jia}
	Let $u$ and $v$ be two positive real numbers such that $\frac{15}{16}\leq u< v<1$. If
	\[u+\sqrt{(1-u)^3}>v,\]
	then
	\[4(1-u)(1-v)-(1-uv)^2>0.\]
\end{lem}

\noindent{\it Proof of Conjecture \ref{cp-1}. }
We first  prove  that $q(n)$ is log-concave for $n\geq 33$.  It is equivalent to proving that for $n\geq 33$,
\[Q(n)\leq1.\]
By Theorem \ref{thm-B}, we see that for  $\nu(n)\geq 67$,
\begin{equation}\label{lem-B-2}
	Q(n)<1-\frac{\pi ^4}{36 \nu(n)^3}+\frac{\pi ^4}{12 \nu(n)^4}-\frac{\pi ^4}{32 \nu(n)^5}+\frac{126+\frac{\pi ^8}{1296} }{\nu(n)^6}.
\end{equation}
It is easy to check that for $\nu(n)\geq 44 $
\[-\frac{\pi ^4}{36 \nu(n)^3}+\frac{\pi ^4}{12 \nu(n)^4}\leq0 \]
and
\[-\frac{\pi ^4}{32 \nu(n)^5}+\frac{126+\frac{\pi ^8}{1296} }{\nu(n)^6}\leq0,\]
we therefore get that $Q(n)<1$ for $n \geq 1365$. It can be verified that $Q(n)<1$ for $ 33\leq  n \leq1365$. Thus, we derive that $Q(n)<1$ for $n\geq 33$, and so  $q(n)$ is log-concave for $n\geq 33$.

To prove that $q(n)$ satisfies the higher order Tur\'an inequalities for $n\geq121$, it suffices to show that for $n\geq 121$,
\begin{equation}\label{eq-q-turan-pf}
4(1-Q(n))(1-Q(n+1))-(1-Q(n)Q(n+1))^2>0.
	\end{equation}
It  can be directly checked that \eqref{eq-q-turan-pf} is true  when $121\leq n\leq 1365$,  so it's enough  to prove that \eqref{eq-q-turan-pf} holds for $n\geq 1365$.
Using Lemma \ref{lem-Jia} as well as the assertion we have proved before that $Q(n+1)<1$ for $n\geq 32$,   we just need to prove that for $n\geq 1365$,
\begin{equation}\label{eq-lem-1}
\frac{15}{16} \leq Q(n)< Q(n+1)
\end{equation}
and
\begin{equation}\label{eq-lem-2}
	Q(n+1)<Q(n)+\sqrt{(1-Q(n))^3} .
\end{equation}

Using Theorem  \ref{thm-B}, we see that for $\nu(n)\geq 67$,
\begin{align}
		Q(n)>1-\frac{\pi ^4}{36 \nu(n)^3}+\frac{\pi ^4}{12 \nu(n)^4}-\frac{\pi ^4}{32 \nu(n)^5}-\frac{135 }{\nu(n)^6}.\nonumber
\end{align}
It's easy to check that for $\nu(n)\geq 5$,
\[\frac{\pi ^4}{12 \nu(n)^4}-\frac{\pi ^4}{32 \nu(n)^5}-\frac{135 }{\nu(n)^6}>0\]
and
\[1-\frac{\pi^4}{36\nu(n)^3}\geq 1-\frac{\pi^4}{36\cdot5^3}>\frac{15}{16}.\]
Hence we obtain that for $\nu(n)\geq 67$,
\[Q(n)>\frac{15}{16}.\]
Using Theorem \ref{thm-B} again, we find that for  $\nu(n)\geq 67$,
	\begin{align}\label{lem-B-1}
	Q(n+1)-Q(n)&>\left(1-\frac{\pi ^4}{36 \nu(n+1)^3}+\frac{\pi ^4}{12 \nu(n+1)^4}-\frac{\pi ^4}{32 \nu(n+1)^5}-\frac{135 }{\nu(n+1)^6}\right)\nonumber\\[6pt]
	&\quad-\left( 1-\frac{\pi ^4}{36 \nu(n)^3}+\frac{\pi ^4}{12 \nu(n)^4}-\frac{\pi ^4}{32 \nu(n)^5}+\frac{126+\frac{\pi ^8}{1296} }{\nu(n)^6}\right).
\end{align}
Note that for $\nu(n)\geq 3$,
\begin{align}
	\left\{
	\begin{aligned}\label{eq-x}
		&\frac{1}{\nu(n+1)^3}<\frac{1}{\nu(n)^3}-\frac{\pi^2}{4\nu(n)^5},\\[6pt]
		&\frac{1}{\nu(n+1)^4}>\frac{1}{\nu(n)^4}-\frac{2\pi^2}{3\nu(n)^6},\\[6pt]
		&\frac{1}{\nu(n+1)^5}<\frac{1}{\nu(n)^5},\\[6pt]
		&	\frac{1}{\nu(n+1)^6}<\frac{1}{\nu(n)^6}.
	\end{aligned}\right.
\end{align}
Applying \eqref{eq-x} to \eqref{lem-B-1}, we obtain that for $\nu(n)\ge 67$,
\begin{align}
	Q(n+1)-Q(n)&>\left(1-\frac{\pi^4}{36\nu(n)^3}+\frac{\pi^4}{12\nu(n)^4}+\frac{-\frac{\pi^4}{32}+\frac{\pi^6}{144}}{\nu(n)^5}-\frac{\frac{\pi^6}{18}+135}{\nu(n)^6}\right)\nonumber\\[6pt]
	&\quad-\left(1-\frac{\pi ^4}{36 \nu(n)^3}+\frac{\pi ^4}{12 \nu(n)^4}-\frac{\pi ^4}{32 \nu(n)^5}+\frac{126+\frac{\pi ^8}{1296} }{\nu(n)^6}\right)\nonumber\\[6pt]
	&=\frac{\pi^6}{144\nu(n)^5}-\frac{261+\frac{\pi ^6}{18}+\frac{\pi ^8}{1296}}{\nu(n)^6}.\nonumber
\end{align}
It is not difficult to check that for $\nu(n)\geq49$,
\[\frac{\pi^6}{144\nu(n)^5}-\frac{261+\frac{\pi ^6}{18}+\frac{\pi ^8}{1296}}{\nu(n)^6}>0,\]
so we get that for $\nu(n)\geq 67$,
\[Q(n+1)-Q(n)>0,\]
and \eqref{eq-lem-1} is verified since $\nu(n)=67$ whenever $n=1365$.

To prove \eqref{eq-lem-2}, using Theorem \ref{thm-B} again, we obtain that for $\nu(n)\geq 67$,
\begin{align}\label{eq-Yn+1-Yn-1}
	Q(n+1)-Q(n)&<\left(	1-\frac{\pi ^4}{36 \nu(n+1)^3}+\frac{\pi ^4}{12 \nu(n+1)^4}-\frac{\pi ^4}{32 \nu(n+1)^5}+\frac{126+\frac{\pi ^8}{1296} }{\nu(n+1)^6}\right)\nonumber\\[6pt]
	&\quad-\left(1-\frac{\pi ^4}{36 \nu(n)^3}+\frac{\pi ^4}{12 \nu(n)^4}-\frac{\pi ^4}{32 \nu(n)^5}-\frac{135 }{\nu(n)^6}\right)\nonumber\\[6pt]
	&=\frac{\pi^4}{36}\left(\frac{1}{\nu(n)^3}-\frac{1}{\nu(n+1)^3}\right)-\frac{\pi^4}{12}\left(\frac{1}{\nu(n)^4}-\frac{1}{\nu(n+1)^4}\right)\nonumber\\[6pt]
	&\quad+\frac{\pi^4}{32}\left(\frac{1}{\nu(n)^5}-\frac{1}{\nu(n+1)^5}\right)+\frac{126+\frac{\pi ^8}{1296}}{\nu(n+1)^6}+\frac{135}{\nu(n)^6}.
\end{align}
It can be checked that for $\nu(n)>0$,
\begin{align}
	\left\{
	\begin{aligned}\label{eq-x-2}
		&\frac{1}{\nu(n)^3}-\frac{1}{\nu(n+1)^3}<\frac{\pi^2}{2\nu(n)^5},\\[6pt]
		&\frac{1}{\nu(n+1)^4}-\frac{1}{\nu(n)^4}<0,\\[6pt]
		&\frac{1}{\nu(n)^5}-\frac{1}{\nu(n+1)^5}<\frac{1}{\nu(n)^5},
	\end{aligned}\right.
\end{align}
and for $\nu(n)\geq 55$,
\begin{align}\label{eq-x-3}
	\frac{126+\frac{\pi ^8}{1296}}{\nu(n+1)^6}+\frac{135}{\nu(n)^6}<\frac{126+\frac{\pi ^8}{1296}+135}{\nu(n)^6}<\frac{\pi^2}{2\nu(n)^5}.
\end{align}
Applying  \eqref{eq-x-2} and \eqref{eq-x-3} to \eqref{eq-Yn+1-Yn-1}, we derive that for $\nu(n)\geq 67$,
\begin{align}\label{eq-Yn+1-Yn-2}
	Q(n+1)-Q(n)&<\frac{\pi^4}{36}\cdot\frac{\pi^2}{2\nu(n)^5}+\frac{\pi^4}{32}\cdot\frac{1}{\nu(n)^5}+\frac{\pi^2}{2\nu(n)^5}\nonumber\\[6pt]
	&=\frac{\frac{\pi^6}{72}+\frac{\pi^4}{32}+\frac{\pi^2}{2}}{\nu(n)^5}.
\end{align}
It remains to show that for $\nu(n)\geq 67$,
\begin{equation}\label{eq-sqrt-1-aa}
	\sqrt{(1-Q(n))^3}>\frac{\frac{\pi^6}{72}+\frac{\pi^4}{32}+\frac{\pi^2}{2}}{\nu(n)^5}.
\end{equation}
Since for  $\nu(n)\geq67$,
\begin{align*}
	1-Q(n)&>\frac{\pi ^4}{36 \nu(n)^3}-\frac{\pi ^4}{12 \nu(n)^4}+\frac{\pi ^4}{32 \nu(n)^5}-\frac{126+\frac{\pi ^8}{1296} }{\nu(n)^6},
\end{align*}
and for $\nu(n)\geq 44$, it can be checked that
\[\frac{\pi ^4}{32 \nu(n)^5}-\frac{126+\frac{\pi ^8}{1296} }{\nu(n)^6}>0\]
and
\[-\frac{\pi ^4}{12 \nu(n)^4}>-\frac{1}{4\nu(n)^3},\]
thus we deduce  that for $\nu(n)\geq 67$,
\begin{equation}\label{eq-sqrt-1}
1-Q(n)>\frac{\pi^4}{36\nu(n)^3}-\frac{1}{4\nu(n)^3}=\frac{\pi^4-9}{36\nu(n)^3}>0.
\end{equation}
It can be easily checked that for $\nu(n)\geq 31$,
\begin{equation}\label{eq-sqrt-2}
	\frac{\sqrt{(\pi^4-9)^3}}{216\nu(n)^{\frac{9}{2}}}>\frac{\frac{\pi^6}{72}+\frac{\pi^4}{32}+\frac{\pi^2}{2}}{\nu(n)^5}.
\end{equation}
Combining   \eqref{eq-sqrt-1} and \eqref{eq-sqrt-2}, we obtain \eqref{eq-sqrt-1-aa}, and so \eqref{eq-lem-2} is valid when $\nu(n)\geq 67$, or equivalently, $n\geq 1365$.  This completes the proof.
\qed

\section{Concluding Remarks}
Let $p_k(n)$ denote the number of partitions of $n$ in which no parts are multiples of $k$. By definition, we see that
\[\sum_{n\geq 0}p_k(n)q^n=\prod_{n=1}^\infty\frac{1-q^{kn}}{1-q^n}.
\]
When $k=2$, this partition function $p_k(n)$ reduces to $q(n)$. In \cite{Craig-Pun-2021},  Craig and Pun also conjectured that the minimal number $N_k$ and $M_k$ such that $p_k(n)$  is log-concave for $n\geq N_k$ and satisfies the higher order Tur\'an inequalities for $n\geq M_k$, where $k=3,4,5$. More precisely, Craig and Pun  \cite{Craig-Pun-2021} conjectured that
\[N_3=58,\quad N_4=17, \quad N_5=42\]
and
\[M_3=185,\quad M_4=64, \quad M_5=137.\]
Based on Chern's asymptotic formulas  for $\eta$-quotients, and using the similar argument in this paper,  we could show that their conjectured values are true. Here we omit the details.  However, our method can not be applied to find the values $N_k$ and $M_k$ for any fixed $k$ such that   $p_k(n)$  is log-concave for $n\geq N_k$ and satisfies the higher order Tur\'an inequalities  for $n\geq M_k$. It would be interesting to find a unified way to determine such $N_k$ and $M_k$ in terms of $k$.

Last but not least, we would like to mention that while studying on higher order Tur\'an inequalities for $p(n)$, Chen \cite{Chen-2017} undertook
a comprehensive study on inequalities pertaining to invariants of a binary form. In particular, he considered the following three invariants of the quartic binary form
\begin{align*}
& A\left(a_0, a_1, a_2, a_3, a_4\right)=a_0 a_4-4 a_1 a_3+3 a_2^2, \\[5pt]
& B\left(a_0, a_1, a_2, a_3, a_4\right)=-a_0 a_2 a_4+a_2^3+a_0 a_3^2+a_1^2 a_4-2 a_1 a_2 a_3, \\[5pt]
&I\left(a_0, a_1, a_2, a_3, a_4\right)= A\left(a_0, a_1, a_2, a_3, a_4\right)^3-27B\left(a_0, a_1, a_2, a_3, a_4\right)^2.
\end{align*}
Chen \cite{Chen-2017} conjectured both the partition function $p(n)$ and the spt-function $\mathrm{spt}(n)$ satisfy the inequalities derived from the  invariants of the quartic binary form  for   large $n$. For the definition of the spt-function, please see Andrews \cite{Andrews-2008} or Chen \cite{Chen-2017}.  Chen's conjectured inequalities  on the partition function  $p(n)$ have recently been proved by Banerjee \cite{Banerjee-2022-1} and Jia and Wang \cite{Jia-Wang-2020}. More precisely, independent work by Banerjee \cite{Banerjee-2022-1} and Jia and Wang \cite{Jia-Wang-2020} proved that $B\left(p({n-1}),p(n), p({n+1}), p({n+2}), p({n+3})\right)>0$ for $n\geq 221$ and Banerjee \cite{Banerjee-2022-1}  showed that $A\left(p({n-1}),p(n), p({n+1}), p({n+2}),\break\right.$
$\left.p({n+3})\right)>0$ for $n\geq185$.

In the same vein, we will present corresponding conjectures on $q(n)$.
\begin{con} Let $a_n=q(n)$, then
\begin{align*}
&A\left(a_{n-1}, a_n, a_{n+1}, a_{n+2}, a_{n+3}\right)>0, & \text { for } n \geq 230, \\[5pt]
&B\left(a_{n-1}, a_n, a_{n+1}, a_{n+2}, a_{n+3}\right)>0, & \text { for } n \geq 272,\\[5pt]
&I\left(a_{n-1}, a_n, a_{n+1}, a_{n+2}, a_{n+3}\right)>0, & \text { for } n \geq 267.
\end{align*}
\end{con}

 \vskip 0.2cm
\noindent{\bf Acknowledgment.} This work
was supported by the National Science Foundation of China.

\end{document}